\documentclass[11pt]{article}
\usepackage{amsmath,amssymb,amsthm}
\usepackage{fancybox}
\usepackage{graphicx}
\usepackage{color}
\usepackage{array}
\usepackage{longtable}
\usepackage{colortab}
\usepackage{colortbl}
\usepackage{arydshln}
\usepackage{tikz}
\usepackage[all]{xy} 

\allowdisplaybreaks     

\begin{document}

\def\fl#1{\left\lfloor#1\right\rfloor}
\def\cl#1{\left\lceil#1\right\rceil}
\def\ang#1{\left\langle#1\right\rangle}
\def\stf#1#2{\left[#1\atop#2\right]}
\def\sts#1#2{\left\{#1\atop#2\right\}}
\def\eul#1#2{\left\langle#1\atop#2\right\rangle}
\def\N{\mathbb N}
\def\Z{\mathbb Z}
\def\R{\mathbb R}
\def\C{\mathbb C}
\def\x{\mathbf x}
\def\y{\mathbf y}
\def\z{\mathbf z}
\newcommand{\ctext}[1]{\raise0.2ex\hbox{\textcircled{\scriptsize{#1}}}}

\newtheorem{theorem}{Theorem}
\newtheorem{Prop}{Proposition}
\newtheorem{Cor}{Corollary}
\newtheorem{Lem}{Lemma}
\newtheorem{Sublem}{Sublemma}
\newtheorem{Def}{Definition}
\newtheorem{Conj}{Conjecture}

\newenvironment{Rem}{\begin{trivlist} \item[\hskip \labelsep{\it
Remark.}]\setlength{\parindent}{0pt}}{\end{trivlist}}

\title{Frobenius numbers associated with Diophantine triples of $x^2+y^2=z^r$ (extended version)
}

\author{
Takao Komatsu
\\
\small Department of Mathematical Sciences, School of Science\\[-0.8ex]
\small Zhejiang Sci-Tech University\\[-0.8ex]
\small Hangzhou 310018 China\\[-0.8ex]
\small \texttt{komatsu@zstu.edu.cn}\\\\ 
Neha Gupta  
\\ 
\small Department of Mathematics, School of Natural Sciences\\[-0.8ex]
\small Shiv Nadar Institute of Eminence\\[-0.8ex]
\small Gautam Buddha Nagar - 201314 India\\[-0.8ex]
\small \texttt{neha.gupta@snu.edu.in}\\\\   
Manoj Upreti
\\ 
\small Department of Mathematics, School of Natural Sciences\\[-0.8ex]
\small Shiv Nadar Institute of Eminence\\[-0.8ex]
\small Gautam Buddha Nagar - 201314 India\\[-0.8ex]
\small \texttt{mu506@snu.edu.in}
}

\date{
\small MR Subject Classifications: Primary 11D07, 11D25; Secondary 05A15,11D04, 20M14
}

\maketitle

\begin{abstract}
We give an explicit formula for the $p$-Frobenius number of triples associated with Diophantine equations $x^2+y^2=z^r$, that is, the largest positive integer that can only be represented in $p$ ways by combining the three integers of the solutions of Diophantine equations $x^2+y^2=z^r$. When $r=2$, the Frobenius number has already been given. 
\\  
{\bf Keywords:} Frobenius problem, Diophantine equations, Pythagorean triples, Ap\'ery set, 
\end{abstract}

\section{Introduction}

Diophantine equations are a fundamental part and one of the oldest branches of number theory. The main study is of polynomial equations or systems of equations, in particular, in integers. Though there are many aspects and applications, 
Diophantine equations are used to characterize certain problems in Diophantine approximations. In \cite{EKS06,EKS07}, we computed upper and lower bounds for the approximation of hyperbolic functions at points $1/s$ ($s=1,2,\dots$) by rationals $x/y$, such that $x$, $y$ and $z$ form Pythagorean triples. In \cite{EKS09,CKL13}, we considered Diophantine approximations $x/y$ to values $\xi$ of hyperbolic functions, where $(x,y,z)$ is the solution of more Diophantine equations, including $x^2+y^2=z^4$.  

For integer $k\ge 2$, consider a set of positive integers $A=\{a_1,\dots,a_k\}$ with $\gcd(A)=\gcd(a_1,\dots,a_k)=1$. To find the number of non-negative integral representations $x_1,x_2,\dots,x_k$, denoted by $d(n;A)=d(n;a_1,a_2,\dots,a_k)$, to $a_1 x_1+a_2 x_2+\dots+a_k x_k=n$ for a given positive integer $n$ is one of the most important and interesting topics. This number is often called the {\it denumerant} and is equal to the coefficient of $x^n$ in $1/(1-x^{a_1})(1-x^{a_2})\cdots(1-x^{a_k})$ (\cite{sy1882}). Sylvester \cite{sy1857} and Cayley \cite{Cayley} showed that $d(n;a_1,a_2,\dots,a_k)$ can be expressed as the sum of a polynomial in $n$ of degree $k-1$ and a periodic function of period $a_1 a_2\cdots a_k$. For two variables, a formula for $d(n;a_1,a_2)$ is obtained in \cite{tr00}. For three variables in the pairwise coprime case $d(n;a_1,a_2,a_3)$, in \cite{Ko03}, the periodic function part is expressed in terms of trigonometric functions.

For a non-negative integer $p$, define $S_p$ and $G_p$ by
$$
S_p(A)=\{n\in\mathbb N_0|d(n;A)>p\}
$$
and
$$
G_p(A)=\{n\in\mathbb N_0|d(n;A)\le p\}
$$
respectively, satisfying $S_p\cup G_p=\mathbb N_0$, which is the set of non-negative integers.
The set $S_p$ is called {\it $p$-numerical semigroup} because $S(A)=S_0(A)$ is a usual numerical semigroup. $G_p$ is the set of {\it $p$-gaps}. Define $g_p(A)$, $n_p(A)$ and $s_p(A)$ by
$$
g_p(A)=\max_{n\in G_p(A)}n,\quad
n_p(A)=\sum_{n\in G_p(A)}1,\quad
s_p(A)=\sum_{n\in G_p(A)}n\,,
$$
respectively, and are called the {\it $p$-Frobenius number}, the {\it $p$-Sylvester number} (or {\it $p$-genus}) and the {\it $p$-Sylvester sum}, respectively.
When $p=0$, $g(A)=g_0(A)$, $n(A)=n_0(A)$ and $s(A)=s_0(A)$ are the original Frobenius number, Sylvester number (or genus) and Sylvester sum, respectively.
Finding such values is one of the crucial matters in the Diophantine problem of Frobenius. More detailed descriptions of the $p$-numerical semigroups and their symmetric properties can be found in \cite{KY24}.

The Frobenius problem (also known as the Coin Exchange Problem or Postage Stamp Problem or Chicken McNugget Problem) has a long history and is one of the popular problems that has attracted the attention of experts as well as amateurs. For two variables $A=\{a,b\}$, it is known that
$$
g(a,b)=(a-1)(b-1)-1\quad\hbox{and}\quad n(a,b)=\frac{(a-1)(b-1)}{2}
$$
(\cite{sy1882,sy1884}).
For three or more variables, the Frobenius number cannot be given by any set of closed formulas, which can be reduced to a finite set of certain polynomials (\cite{cu90}).
For three variables, various algorithms have been devised for finding the Frobenius number. For example, in \cite{RGS04}, the Frobenius number is uniquely determined by six positive integers that are the solution to a system of three polynomial equations. In \cite{Fel06}, a general algorithm is given by using a $3\times 3$ matrix.
Nevertheless, explicit closed formulas have been found only for some special cases, including arithmetic, geometric, Mersenne, repunits and triangular (see \cite{RR18,RBT2015,RBT2017} and references therein). We are interested in finding explicit closed forms, which is one of the most crucial matters in Frobenius problem. Our method has an advantage in terms of visually grasping the elements of the Ap\'ery set, and is more useful to get more related values, including genus (Sylvester number), Sylvester sum \cite{tr08}, weighted power Sylvester sum \cite{Ko22,KZ0,KZ} and so on.

We are interested in finding a closed or explicit form for the $p$-Frobenius number, which is more difficult when $p>0$. For three or more variables, no concrete examples had been found until recently, when the first author succeeded in giving the $p$-Frobenius number as a closed-form expression for the triangular number triplet (\cite{Ko22a}), for repunits (\cite{Ko22b}), Fibonacci triplet (\cite{Ko23b}), Jacobsthal triplets (\cite{KP,KLP}) and arithmetic triplets (\cite{KY23c}).

We give an explicit formula for the $p$-Frobenius number of triples associated with Diophantine equations, that is the largest positive integer that can only be represented in $p$ ways by combining the three integers of the solutions of Diophantine equations. 
When $p=0$, it is the original Frobenius number in the famous Diophantine problem of Frobenius. We also obtain closed forms for the number of positive integers, and the largest positive integer that can be represented in only $p$ ways by combining the three integers of the Diophantine triple. 

In this paper, we study the numerical semigroup of the triples $(x,y,z)$, satisfying the Diophantine equation $x^2+y^2=z^r$ ($r\ge 2$). When $r=2$, the Frobenius number of the Pythagorean triple is given in \cite{GHKKLLNPP}. The Frobenius number of a little modified triple is studied in \cite{KS}.


\section{Preliminaries}

We introduce the $p$-Ap\'ery set (see \cite{Apery}) below in order to obtain the formulas for $g_p(A)$, $n_p(A)$ and $s_p(A)$. Without loss of generality, we assume that $a_1=\min(A)$.

\begin{Def}
Let $p$ be a nonnegative integer. For a set of positive integers $A=\{a_1,a_2,\dots,a_\kappa\}$ with $\gcd(A)=1$ and $a_1=\min(A)$ we denote by
$$
{\rm Ap}_p(A)={\rm Ap}_p(a_1,a_2,\dots,a_\kappa)=\{m_0^{(p)},m_1^{(p)},\dots,m_{a_1-1}^{(p)}\}\,,
$$
the $p$-Ap\'ery set of $A$, where each positive integer $m_i^{(p)}$ $(0\le i\le a_1-1)$ satisfies the conditions:
$$
{\rm (i)}\, m_i^{(p)}\equiv i\pmod{a_1},\quad{\rm (ii)}\, m_i^{(p)}\in S_p(A),\quad{\rm (iii)}\, m_i^{(p)}-a_1\not\in S_p(A)
$$
Note that $m_0^{(0)}$ is defined to be $0$.
\label{apery}
\end{Def}

\noindent
It follows that for each $p$,
$$
{\rm Ap}_p(A)\equiv\{0,1,\dots,a_1-1\}\pmod{a_1}\,.
$$

When $k\ge 3$, it is hard to find any explicit form of $g_p(A)$ as well as $n_p(A)$ and $s_p(A)$.  Nevertheless, the following convenient formulas are known (for a more general case, see \cite{Ko-p}). Though finding $m_j^{(p)}$ is enough hard in general, we can obtain it for some special sequences $(a_1,a_2,\dots,a_k)$. 

\begin{Lem}
Let $k$ and $p$ be integers with $k\ge 2$ and $p\ge 0$.
Assume that $\gcd(a_1,a_2,\dots,a_k)=1$.  We have
\begin{align}
g_p(a_1,a_2,\dots,a_k)&=\left(\max_{0\le j\le a_1-1}m_j^{(p)}\right)-a_1\,,
\label{mp-g}
\\
n_p(a_1,a_2,\dots,a_k)&=\frac{1}{a_1}\sum_{j=0}^{a_1-1}m_j^{(p)}-\frac{a_1-1}{2}\,,
\label{mp-n}
\\
s_p(a_1,a_2,\dots,a_k)&=\frac{1}{2 a_1}\sum_{j=0}^{a_1-1}\bigl(m_j^{(p)}\bigr)^2-\frac{1}{2}\sum_{j=0}^{a_1-1}m_j^{(p)}+\frac{a_1^2-1}{12}\,.
\label{mp-s}
\end{align}
\label{lem-mp}
\end{Lem}

\noindent
{\it Remark.}
When $p=0$, the formulas (\ref{mp-g}), (\ref{mp-n}) and (\ref{mp-s}) reduce to the formulas by Brauer and Shockley \cite{bs62}, Selmer \cite{se77}, and Tripathi \cite{tr08}, respectively:
\begin{align*}
g(a_1,a_2,\dots,a_k)&=\left(\max_{1\le j\le a_1-1}m_j\right)-a_1\,,\\
n(a_1,a_2,\dots,a_k)&=\frac{1}{a_1}\sum_{j=1}^{a_1-1}m_j-\frac{a_1-1}{2}\,,\\
s(a_1,a_2,\dots,a_k)&=\frac{1}{2 a_1}\sum_{j=1}^{a_1-1}(m_j)^2-\frac{1}{2}\sum_{j=1}^{a_1-1}m_j+\frac{a_1^2-1}{12}\,,
\end{align*}
where $m_j=m_j^{(0)}$ ($1\le j\le a_1-1$) with $m_0=m_0^{(0)}=0$.
More general formulas using Bernoulli numbers can be seen in \cite{Ko22}.

\section{$x^2+y^2=z^3$} 

For the solution of the Diophantine equation $x^2+y^2=z^r$, we obtain the parameterization 
\begin{align*}
x&=\sum_{k=0}^{\fl{r/2}}(-1)^k\binom{r}{2 k}s^{r-2 k}t^{2 k}\,,\\
y&=\sum_{k=0}^{\fl{(r-1)/2}}(-1)^k\binom{r}{2 k+1}s^{r-2 k-1}t^{2 k+1}\,,\\
z&=s^2+t^2\,,  
\end{align*}
where $s$ and $t$ are of opposite parity with $\gcd(s,t)=1$.  

The case where $r=2$ has already been discussed in \cite{GHKKLLNPP,KS}. 
Namely, 
$$
g_0(s^2-t^2,2 s t,s^2+t^2)=(s-1)(s^2-t^2)+(s-1)(2 s t)-(s^2+t^2)\,. 
$$

Let $r=3$. The triple of the Diophantine equation $x^2+y^2=z^3$ is parameterized as  
$$
(x,y,z)=\bigl(s(s^2-3 t^2),t(3 s^2-t^2),s^2+t^2\bigr)\,. 
$$ 
For convenience, we put 
\begin{equation}
\x:=s(s^2-3 t^2),\quad \y:=t(3 s^2-t^2),\quad \z:=s^2+t^2\,. 
\label{xyz_r=3}
\end{equation} 
Since $\x,\y,\z>0$ and $\gcd(\x,\y,\z)=1$, we see that $s>\sqrt{3}t$, $\gcd(s,t)=1$ and $s\not\equiv t\pmod{2}$.

When $\x>\z$, the Frobenius number of this triple is given as follows.   

\begin{theorem}  
\begin{align*}
&g_0\bigl(s(s^2-3 t^2),t(3 s^2-t^2),s^2+t^2\bigr)\\
&=\begin{cases}
&(s-1)s(s^2-3 t^2)+(s-1)t(3 s^2-t^2)-(s^2+t^2)\\
&\qquad\qquad\qquad\text{if $s>(1+\sqrt{2})t$};\\ 
&(2 s+t-1)s(s^2-3 t^2)+(t-1)t(3 s^2-t^2)-(s^2+t^2)\\
&\qquad\qquad\qquad\text{if $2 t<s<(1+\sqrt{2})t$};\\ 
&(2 s+t-1)s(s^2-3 t^2)+(s-t-1)t(3 s^2-t^2)-(s^2+t^2)\\
&\qquad\qquad\qquad\text{if $(2+\sqrt{13})t/3<s<2 t$};\\ 
&(5 s+3 t-1)s(s^2-3 t^2)+(2 t-s-1)t(3 s^2-t^2)-(s^2+t^2)\\
&\qquad\qquad\qquad\text{if $1.8139 t<s<(2+\sqrt{13})t/3$};\\ 
&(2 s+t-1)s(s^2-3 t^2)+(2 s-3 t-1)t(3 s^2-t^2)-(s^2+t^2)\\
&\qquad\qquad\qquad\text{if $(3+\sqrt{34})t/5<s<1.8139 t$};\\ 
&(7 s+4 t-1)s(s^2-3 t^2)+(2 t-s-1)t(3 s^2-t^2)-(s^2+t^2)\\
&\qquad\qquad\qquad\text{if $\sqrt{3}\,t<s<(3+\sqrt{34})t/5$}\,.  
\end{cases}
\end{align*}
Here, $1.8139\dots$ is the positive real root of $3 x^4-7 x^3+6 x^2-3 x-5=0$. 
\label{th:2-2-3-0}
\end{theorem}  

\noindent 
{\it Remark.}  
When $\x<\z$, that is, $s\sqrt{(s-1)/(3 s+1)}<t<s/\sqrt{3}$, there is no uniform pattern for the Frobenius number. We need a separate discussion for each case. See below for details.

\subsection{The case where $\sqrt{3}t<s<(2+\sqrt{3})t$}  

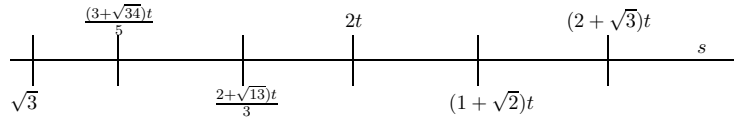
\begin{table}[htbp]
  \centering
\scalebox{0.7}{
\begin{tabular}{ccccccccccccc}
&&&\!\!\!\!\!\!\!\!\!\!\!\!\!\!\!\!\!\!$\frac{(3+\sqrt{34})t}{5}$&&&&\!\!\!\!\!\!\!\!\!\!\!\!\!\!\!\!\!\!$2 t$&&&&\!\!\!\!\!\!\!\!\!\!\!\!\!\!\!\!\!\!$(2+\sqrt{3})t$&\\ 
\multicolumn{1}{c|}{}&&\multicolumn{1}{c|}{}&&\multicolumn{1}{c|}{}&&\multicolumn{1}{c|}{}&&\multicolumn{1}{c|}{}&&\multicolumn{1}{c|}{}&&$s$\\ 
\cline{1-2}\cline{3-4}\cline{5-6}\cline{7-8}\cline{9-10}\cline{11-12}\cline{13-13}
\multicolumn{1}{c|}{}&&\multicolumn{1}{c|}{}&&\multicolumn{1}{c|}{}&&\multicolumn{1}{c|}{}&&\multicolumn{1}{c|}{}&&\multicolumn{1}{c|}{}&&\\ 
&\!\!\!\!\!\!\!\!\!\!\!\!$\sqrt{3}$&\phantom{\qquad}&\phantom{\qquad}&\phantom{\qquad}&\!\!\!\!\!\!\!\!\!\!\!\!$\frac{2+\sqrt{13})t}{3}$&\phantom{\qquad}&\phantom{\qquad}&\phantom{\qquad}&\!\!\!\!\!\!\!\!\!\!\!\!$(1+\sqrt{2})t$&\phantom{\qquad}&\phantom{\qquad}&\phantom{\qquad}\\
\end{tabular}
}
\end{table}

If $\sqrt{3}t<s<(2+\sqrt{3})t$, then $0<\x<\y$. Hence, $\x<\z<\y$ or $\z<\x<\y$.

First, let $\z<\x<\y$. Since $(2 s+t)\x+(2 t-s)\y=2(s^2-s t-t^2)\z$ with $s^2-s t-t^2>(2-\sqrt{3})t^2$, 
\begin{equation}
(2 s+t)\x+(2 t-s)\y\equiv\z\quad\hbox{and}\quad (2 s+t)\x+(2 t-s)\y>0\,. 
\label{eq:111}
\end{equation} 

\noindent 
$\bullet$ The case where $(1+\sqrt{2})t<s<(2+\sqrt{3})t$\\ 
\noindent 
Then 
the elements of the ($0$-)Ap\'ery set are given as in Table \ref{tb:223-3}, where each point $(X,Y)$ corresponds to the expression $X\x+Y\y$ and the area of the ($0$-)Ap\'ery set is equal to $\z=s^2+t^2$.  

\begin{table}[htbp]
  \centering
\scalebox{0.7}{
\begin{tabular}{ccccccc}
\cline{1-2}\cline{3-4}\cline{5-6}\cline{6-7}
\multicolumn{1}{|c}{$(0,0)$}&$\cdots$&$(s-1,0)$&$(s,0)$&$\cdots$&$\cdots$&\multicolumn{1}{c|}{$(s+t-1,0)$}\\
\multicolumn{1}{|c}{$\vdots$}&&$\vdots$&$\vdots$&&&\multicolumn{1}{c|}{$\vdots$}\\
\multicolumn{1}{|c}{$(0,t-1)$}&$\cdots$&$(s-1,t-1)$&$(s,t-1)$&$\cdots$&$\cdots$&\multicolumn{1}{c|}{$(s+t-1,t-1)$}\\
\cline{4-5}\cline{6-7}
\multicolumn{1}{|c}{$(0,t)$}&$\cdots$&\multicolumn{1}{c|}{$(s-1,t)$}&&&&\\\multicolumn{1}{|c}{$\vdots$}&&\multicolumn{1}{c|}{$\vdots$}&&&&\\
\multicolumn{1}{|c}{$\vdots$}&&\multicolumn{1}{c|}{$\vdots$}&&&&\\
\multicolumn{1}{|c}{$(0,s-1)$}&$\cdots$&\multicolumn{1}{c|}{$(s-1,s-1)$}&&&&\\
\cline{1-2}\cline{3-3}
\end{tabular}
}
  \caption{${\rm Ap}_0(\x ,\y ,\z )$ when $(1+\sqrt{2})t<s<(2+\sqrt{3})t$}
  \label{tb:223-3}
\end{table}

Since $(1+\sqrt{2})t<s<(2+\sqrt{3})t$, we see that $(s+t)\x>(s-t)\y$.  
By $(s+t)\x\equiv(s-t)\y\pmod{\z }$ and by $s\x\equiv -t\y\pmod{\z}$ 
the sequence $\{\ell\x\pmod{\z}\}_{\ell=0}^{\z-1}$ is given by 
\begin{align}
&(0,0),(1,0),\dots,(s+t-1,0),(0,s-t),(1,s-t),\dots,(s-1,s-t),\notag\\
&(0,s-2 t),(1,s-2 t),\dots,(s+t-1,s-2 t),\notag\\
&(0,2 s-3 t),(1,2 s-3 t),\dots,(s-1,2 s-3 t),\notag\\
&(0,2 s-4 t),(1,2 s-4 t),\dots,(s-1,2 s-4 t),\notag\\
&(0,2 s-5 t),(1,2 s-5 t),\dots,(s+t-1,2 s-5 t)\notag\\
&(0,3 s-6 t),(1,3 s-6 t),\dots,(s-1,3 s-6 t),\notag\\
&(0,3 s-7 t),(1,3 s-7 t),\dots,(s-1,3 s-7 t)\notag\\
&\dots,(s-1,s t-(s-1)t)\,.
\label{sq:223-3}
\end{align}
Namely, after $(s-1,s t-(s-1)t)$, the next point adding $\x \pmod{\z }$ is $(0,0)$. 
Note that the typical patterns in the sequence (\ref{sq:223-3}) are shown as follows: if $k_1 s-k_2 t\le t-1$
$$ 
(0,k_1 s-k_2 t),(1,k_1 s-k_2 t),\dots,(s+t-1,k_1 s-k_2 t),(0,(k_1+1)s-(k_2+1)t)
$$
and if $k_1 s-k_2 t\ge t$
$$ 
(0,k_1 s-k_2 t),(1,k_1 s-k_2 t),\dots,(s-1,k_1 s-k_2 t),(0,k_1 s-(k_2+1)t)\,. 
$$
Since $\gcd(s,t)=1$, all the points inside of the area in Table \ref{tb:223-3} appear in the sequence (\ref{sq:223-3}) just once.  
Since $\gcd(\x ,\z )=1$, the sequence $\{\ell\x\pmod{\z }\}_{\ell=0}^{\z -1}$ is equivalent to the sequence $\{\ell\pmod{\z }\}_{\ell=0}^{\z -1}$. 

Comparing the elements at $(s+t-1,t-1)$ and $(s-1,s-1)$, taking possible maximal values, we find that the element at $(s-1,s-1)$ is the largest in the Ap\'ery set because 
$$
(s-1)\x+(s-1)\y-\bigl((s+t-1)\x+(t-1)\y\bigr)=t\bigl(s^2(2 s-3)+(2 s+t)t^2\bigr)>0\,. 
$$ 
By Lemma \ref{lem-mp} (\ref{mp-g}), we have 
\begin{align*}
g\bigl(\x ,\y ,\z \bigr)
&=(s-1)\x +(s-1)\y -\z\\
&=(s-1)(s-t)(s^2+4 s t+t^2)-(s^2+t^2)\,. 
\end{align*}   
\bigskip

\noindent 
$\bullet$ The case where $2 t<s<(\sqrt{2}+1)t$\\
\noindent 
Then $(s+t)\x \equiv(s-t)\y \pmod{\z }$ but $(s+t)\x <(s-t)\y $. 
Nevertheless, by $(2 s+t)\x-(s-2 t)\y=2(s^2-s t-t^2)\z>0$, we have $(2 s+t)\x \equiv(s-2 t)\y \pmod{\z }$ and $(2 s+t)\x >(s-2 t)\y$.  
For example, $(s,t)=(9,4)$ satisfies this condition, so $(x,y,z)=(297,908,97)$. By a similar way, we know that 
all the elements of the ($0$-)Ap\'ery set are given as in Table \ref{tb:223-4}.

\begin{table}[htbp]
  \centering
\scalebox{0.7}{
\begin{tabular}{ccccccc}
\cline{1-2}\cline{3-4}\cline{5-6}\cline{6-7}
\multicolumn{1}{|c}{$(0,0)$}&$\cdots$&$(s-1,0)$&$(s,0)$&$\cdots$&$\cdots$&\multicolumn{1}{c|}{$(2 s+t-1,0)$}\\
\multicolumn{1}{|c}{$\vdots$}&&$\vdots$&$\vdots$&&&\multicolumn{1}{c|}{$\vdots$}\\
\multicolumn{1}{|c}{$(0,t-1)$}&$\cdots$&$(s-1,t-1)$&$(s,t-1)$&$\cdots$&$\cdots$&\multicolumn{1}{c|}{$(2 s+t-1,t-1)$}\\
\cline{4-5}\cline{6-7}
\multicolumn{1}{|c}{$(0,t)$}&$\cdots$&\multicolumn{1}{c|}{$(s-1,t)$}&&&&\\\multicolumn{1}{|c}{$\vdots$}&&\multicolumn{1}{c|}{$\vdots$}&&&&\\
\multicolumn{1}{|c}{$\vdots$}&&\multicolumn{1}{c|}{$\vdots$}&&&&\\
\multicolumn{1}{|c}{$(0,s-t-1)$}&$\cdots$&\multicolumn{1}{c|}{$(s-1,s-t-1)$}&&&&\\
\cline{1-2}\cline{3-3}
\end{tabular}
}
  \caption{${\rm Ap}_0(\x ,\y ,\z )$ when $2 t<s<(\sqrt{2}+1)t$}
  \label{tb:223-4}
\end{table}

Compare the elements at $(2 s+t-1,t-1)$ and $(s-1,s-t-1)$, which take possible maximal values. Since the real roots of $-x^4+2 x^3-3 x^2+2 x+2=0$ are $-0.4909$ and $1.4909$, together with $s>2 t$, we see that
$$
(s-1)\x+(s-t-1)\y-\bigl((2 s+t-1)\x+(t-1)\y\bigr)=-s^4+2 s^3 t 3 s^2 t^2+2 s t^3+2 t^4<0\,, 
$$ 
we find that the element at $(2 s+t-1,t-1)$ is the largest in the Ap\'ery set.  
By Lemma \ref{lem-mp} (\ref{mp-g}), we have 
$$ 
g\bigl(\x ,\y ,\z \bigr)
=(2 s+t-1)\x +(t-1)\y -\z\,. 
$$   
\bigskip

\noindent 
$\bullet$ The case where $(2+\sqrt{13})t/3<s<2 t$\\  
\noindent 
For example, $(s,t)=(27,14)$ satisfies this condition, so $(x,y,z)=(3807,27874,925)$.   

By $(2+\sqrt{13})t/3<s$, we have $(3 s+2 t)\x-(2 s-3 t)\y=(3 s^2-4 s t-3 t^2)\z>0$. So, $(3 s+2 t)\x \equiv(2 s-3 t)\y \pmod{\z }$ and $(3 s+2 t)\x >(2 s-3 t)\y$. Together with (\ref{eq:111}), all the elements of the ($0$-)Ap\'ery set are given as in Table \ref{tb:223-5}. 
  
\begin{table}[htbp]
  \centering
\scalebox{0.7}{
\begin{tabular}{ccccccc}
\cline{1-2}\cline{3-4}\cline{5-6}\cline{6-7}
\multicolumn{1}{|c}{$(0,0)$}&$\cdots$&$(2 s+t-1,0)$&$(2 s+t,0)$&$\cdots$&$\cdots$&\multicolumn{1}{c|}{$(3 s+2 t-1,0)$}\\
\multicolumn{1}{|c}{$\vdots$}&&$\vdots$&$\vdots$&&&\multicolumn{1}{c|}{$\vdots$}\\
\multicolumn{1}{|c}{$(0,2t-s-1)$}&$\cdots$&$(2 s+t-1,2 t-s-1)$&$(2 s+t,2 t-s-1)$&$\cdots$&$\cdots$&\multicolumn{1}{c|}{$(3 s+2 t-1,2 t-s-1)$}\\
\cline{4-5}\cline{6-7}
\multicolumn{1}{|c}{$(0,2 t-s)$}&$\cdots$&\multicolumn{1}{c|}{$(2 s+t-1,2 t-s)$}&&&&\\\multicolumn{1}{|c}{$\vdots$}&&\multicolumn{1}{c|}{$\vdots$}&&&&\\
\multicolumn{1}{|c}{$\vdots$}&&\multicolumn{1}{c|}{$\vdots$}&&&&\\
\multicolumn{1}{|c}{$(0,s-t-1)$}&$\cdots$&\multicolumn{1}{c|}{$(2 s+t-1,s-t-1)$}&&&&\\
\cline{1-2}\cline{3-3}
\end{tabular}
}
  \caption{${\rm Ap}_0(\x ,\y ,\z )$ when $(2+\sqrt{13})t/3<s<2 t$}
  \label{tb:223-5}
\end{table}

Compare the elements at $(3 s+2 t-1,2 t-s-1)$ and $(2 s+t-1,s-t-1)$, which take the possible maximal values.  We find that the element at $(s-1,s-1)$ is the largest in the Ap\'ery set because $1.8685=(2+\sqrt{13})t/3<s<2 t$ and 
\begin{align*}
&(2 s+t-1)\x+(s-t-1)\y-\bigl((3 s+2 t-1)\x+(2 t-s-1)\y\bigr)\\
&=-(s^4- 5 s^3 t+6 s^2 t^2-s t^3-3 t^4)>0 
\end{align*} 
for $-0.5268\, t<s<3.3968\, t$.  
By Lemma \ref{lem-mp} (\ref{mp-g}), we have 
$$
g\bigl(\x ,\y ,\z \bigr)
=(2 s+t-1)\x +(s-t-1)\y -\z\,. 
$$   
\bigskip

\noindent 
$\bullet$ The case where $(3+\sqrt{34})t/5<s<(2+\sqrt{13})t/3$\\  
\noindent 
For example, $(s,t)=(24,13)$ satisfies this condition, so $(x,y,z)=(1656,20267,745)$. 

In this case, $(3 s+2 t)\x \equiv(2 s-3 t)\y \pmod{\z }$ but $(3 s+2 t)\x <(2 s-3 t)\y $. 
Nevertheless, by $(3+\sqrt{34})t/5<s$ we have $(5 s+3 t)\x-(3 s-5 t)\y=(5 s^2-6 s t-5 t^2)\z>0$. So, $(5 s+3 t)\x \equiv(3 s-5 t)\y \pmod{\z }$ and $(5 s+3 t)\x >(3 s-5 t)\y$. 
Together with (\ref{eq:111}), all the elements of the ($0$-)Ap\'ery set are given as in Table \ref{tb:223-6}.  

\begin{table}[htbp]
  \centering
\scalebox{0.7}{
\begin{tabular}{ccccccc}
\cline{1-2}\cline{3-4}\cline{5-6}\cline{6-7}
\multicolumn{1}{|c}{$(0,0)$}&$\cdots$&$(2 s+t-1,0)$&$(2 s+t,0)$&$\cdots$&$\cdots$&\multicolumn{1}{c|}{$(5 s+3 t-1,0)$}\\
\multicolumn{1}{|c}{$\vdots$}&&$\vdots$&$\vdots$&&&\multicolumn{1}{c|}{$\vdots$}\\
\multicolumn{1}{|c}{$(0,2t-s-1)$}&$\cdots$&$(2 s+t-1,2 t-s-1)$&$(2 s+t,2 t-s-1)$&$\cdots$&$\cdots$&\multicolumn{1}{c|}{$(5 s+3 t-1,2 t-s-1)$}\\
\cline{4-5}\cline{6-7}
\multicolumn{1}{|c}{$(0,2 t-s)$}&$\cdots$&\multicolumn{1}{c|}{$(2 s+t-1,2 t-s)$}&&&&\\\multicolumn{1}{|c}{$\vdots$}&&\multicolumn{1}{c|}{$\vdots$}&&&&\\
\multicolumn{1}{|c}{$\vdots$}&&\multicolumn{1}{c|}{$\vdots$}&&&&\\
\multicolumn{1}{|c}{$(0,2 s-3 t-1)$}&$\cdots$&\multicolumn{1}{c|}{$(2 s+t-1,2 s-3 t-1)$}&&&&\\
\cline{1-2}\cline{3-3}
\end{tabular}
}
  \caption{${\rm Ap}_0(\x ,\y ,\z )$ when $(3+\sqrt{34})t/5<s<(2+\sqrt{13})t/3$}
  \label{tb:223-6}
\end{table}

Comparing the elements at $(5 s+3 t-1,2 t-s-1)$ and $(2 s+t-1,2 s-3 t-1)$, taking possible maximal values, we find that there are two possibilities. 
By Lemma \ref{lem-mp} (\ref{mp-g}), if 
\begin{align*}
&(2 s+t-1)\x+(2 s-3 t-1)\y-\bigl((5 s+3 t-1)\x+(2 t-s-1)\y\bigr)\\
&=-3 s^4+7 s^3 t-6 s^2 t^2+3 s t^3+5 t^4>0\,, 
\end{align*}
which is equivalent to $-0.5553\, t<s<1.8139\, t$, namely if $1.7661\, t=(3+\sqrt{34})t/5<s<1.8139\, t$,  
then 
$$ 
g\bigl(\x ,\y ,\z \bigr)
=(2 s+t-1)\x +(2 s-3 t-1)\y -\z\,. 
$$   
Otherwise, namely, if $1.8139\, t<s<(2+\sqrt{13})t/3=1.8685\, t$
$$ 
g\bigl(\x ,\y ,\z \bigr)
=(5 s+3 t-1)\x +(2 t-s-1)\y -\z\,. 
$$   
\bigskip

\noindent 
$\bullet$ The case where $\sqrt{3}t<s<(3+\sqrt{34})t/5$\\   
\noindent 
For example, $(s,t)=(44,25)$ satisfies this condition, so $(x,y,z)=(2684,129575,2561)$.  

In this case, $(5 s+3 t)\x \equiv(3 s-5 t)\y \pmod{\z }$ but $(5 s+3 t)\x <(3 s-5 t)\y $.  
Nevertheless, by $s>\sqrt{3}t=1.732\,t>(4+\sqrt{65})t/7=1.723\,t$, we have $(7 s+4 t)\x-(4 s-7 t)\y=(7 s^2 - 8 s t - 7 t^2)\z>0$. So, $(7 s+4 t)\x \equiv(4 s-7 t)\y \pmod{\z }$ and $(7 s+4 t)\x>(4 s-7 t)\y$.   
Together with (\ref{eq:111}), all the elements of the ($0$-)Ap\'ery set are given as in Table \ref{tb:223-7}.  
 
\begin{table}[htbp]
  \centering
\scalebox{0.7}{
\begin{tabular}{ccccccc}
\cline{1-2}\cline{3-4}\cline{5-6}\cline{6-7}
\multicolumn{1}{|c}{$(0,0)$}&$\cdots$&$(2 s+t-1,0)$&$(2 s+t,0)$&$\cdots$&$\cdots$&\multicolumn{1}{c|}{$(7 s+4 t-1,0)$}\\
\multicolumn{1}{|c}{$\vdots$}&&$\vdots$&$\vdots$&&&\multicolumn{1}{c|}{$\vdots$}\\
\multicolumn{1}{|c}{$(0,2t-s-1)$}&$\cdots$&$(2 s+t-1,2 t-s-1)$&$(2 s+t,2 t-s-1)$&$\cdots$&$\cdots$&\multicolumn{1}{c|}{$(7 s+4 t-1,2 t-s-1)$}\\
\cline{4-5}\cline{6-7}
\multicolumn{1}{|c}{$(0,2 t-s)$}&$\cdots$&\multicolumn{1}{c|}{$(2 s+t-1,2 t-s)$}&&&&\\\multicolumn{1}{|c}{$\vdots$}&&\multicolumn{1}{c|}{$\vdots$}&&&&\\
\multicolumn{1}{|c}{$\vdots$}&&\multicolumn{1}{c|}{$\vdots$}&&&&\\
\multicolumn{1}{|c}{$(0,3 s-5 t-1)$}&$\cdots$&\multicolumn{1}{c|}{$(2 s+t-1,3 s-5 t-1)$}&&&&\\
\cline{1-2}\cline{3-3}
\end{tabular}
}
  \caption{${\rm Ap}_0(\x ,\y ,\z )$ when $\sqrt{3}t<s<(3+\sqrt{34})t/5$}
  \label{tb:223-7}
\end{table}
 
Comparing the elements at $(7 s+4 t-1,2 t-s-1)$ and $(2 s+t-1,3 s-5 t-1)$, taking possible maximal values, we find that the element at $(7 s+4 t-1,2 t-s-1)$ is the largest in the Ap\'ery set because by $\sqrt{3}t<s<(3+\sqrt{34})t/5$ 
\begin{align*}
&(2 s+t-1)\x+(3 s-5 t-1)\y-\bigl((7 s+4 t-1)\x+(2 t-s-1)\y\bigr)\\
&=-5 s^4+9 s^3 t-6 s^2 t^2+5 s t^3+7 t^4<0\,. 
\end{align*}
Note that $-5 x^4+9 x^3-6 x^2+5 x+7=0$ has the real roots at $0.5702$ and $1.71692$ and $s>\sqrt{3}t=1.732 t$.  
By Lemma \ref{lem-mp} (\ref{mp-g}), we have 
$$ 
g\bigl(\x ,\y ,\z \bigr)
=(7 s+4 t-1)\x +(2 t-s-1)\y -\z\,. 
$$    
\bigskip

\subsection{The case where $\x<\z<\y$} 

If $\x<\z<\y$, then by $s(s^2-3 t^2)<\z $ and $s>\sqrt{3}t$, we have 
$$
s\sqrt{\frac{s-1}{3 s+1}}<t<\frac{s}{\sqrt{3}}\,. 
$$
As both $s$ and $t$ are positive integers, we can find 
\begin{align*}  
&(s,t)=(2,1),(7,4),(9,5),(11,6),(14,8),(16,9),(21,12),(23,13),(26,15),\\
&\qquad (28,16),(30,17),(33,19),(35,20),(37,21),(40,23),\dots\,. 
\end{align*}
With additional conditions $\gcd(s,t)=1$ and $s\not\equiv t\pmod{2}$, the suitable pairs are given by 
$$
(s,t)=(2,1),(7,4),(11,6),(16,9),(26,15),(30,17),(40,23),\dots\,. 
$$ 
Notice that the area of the ($0$-)Ap\'ery set is equal to $\x=s(s^2-3 t^2)$. 

In this case, there is no uniform pattern. We have to calculate the Frobenius number for each concrete case. For example, when $(s,t)=(2,1),(7,4),(26,15)$, the area of the Ap\'ery set is given as in Table \ref{tb:223-11}, where each point $(Z,Y)$ corresponds to the expression $Z \z +Y \y $ and the area of the ($0$-)Ap\'ery set is equal to $\x=s(s^2-3 t^2)$.  

\begin{table}[htbp]
  \centering
\begin{tabular}{ccc}
\cline{1-2}\cline{3-3}
\multicolumn{1}{|c}{$(0,0)$}&$\cdots$&\multicolumn{1}{c|}{$(z_0,0)$}\\
\cline{1-2}\cline{3-3}
\end{tabular}
  \caption{${\rm Ap}_0(\x ,\y ,\z )$ when $(s,t)=(2,1),(7,4),(25,15),\dots$}
  \label{tb:223-11}
\end{table}

Here, when $(s,t)=(2,1),(7,4),(26,15)$, $A=2,7,26$, respectively. Then the Frobenius number is given by 
$$
g(\x,\y,\z)=z_0(s^2+t^2)-s(s^2-3 t^2)\,. 
$$

When $(s,t)=(11,6),(16,9),(30,17),(40,23)$, the area of the Ap\'ery set is given as in Table \ref{tb:223-12}. 

\begin{table}[htbp]
  \centering
\scalebox{0.9}{
\begin{tabular}{ccccc}
\cline{1-2}\cline{3-4}\cline{5-5}
\multicolumn{1}{|c}{$(0,0)$}&$\cdots$&$\cdots$&$\cdots$&\multicolumn{1}{c|}{}\\
\multicolumn{1}{|c}{}&$\cdots$&$\cdots$&$\cdots$&\multicolumn{1}{c|}{$(z_1,y_1)$}\\
\cline{4-5}
\multicolumn{1}{|c}{}&$\cdots$&\multicolumn{1}{c|}{$\cdots$}&&\\
\multicolumn{1}{|c}{}&$\cdots$&\multicolumn{1}{c|}{$(z_2,y_2)$}&&\\
\cline{1-2}\cline{3-3}
\end{tabular}
}
  \caption{${\rm Ap}_0(\x ,\y ,\z )$ when $(s,t)=(11,6),(16,9),(30,17),(40,23),\dots$}
  \label{tb:223-12}
\end{table}

Here, 
\begin{align*}
(s,t)=(11,6)&\Longrightarrow~(z_1,y_1)=(37,2),~(z_2,y_2)=(28,3);\\
(s,t)=(16,9)&\Longrightarrow~(z_1,y_1)=(61,2),~(z_2,y_2)=(10,4);\\
(s,t)=(30,17)&\Longrightarrow~(z_1,y_1)=(201,0),~(z_2,y_2)=(196,4);\\
(s,t)=(40,23)&\Longrightarrow~(z_1,y_1)=(136,2),~(z_2,y_2)=(108,3)
\end{align*}
As seen, there seems to be no regularity. So, one has to find the values of $(z_p,y_p)$ one by one.

When $(s,t)=(11,6),(30,17),(40,23)$, the element at $(z_1,y_1)$ is bigger than the element at $(z_2,y_2)$. Then the Frobenius number is given by 
$$
g(\x,\y,z)=z_1(s^2+t^2)+y_1 t(3 s^2-t^2)-s(s^2-3 t^2)\,. 
$$
When $(s,t)=(16,9)$, the element at $(z_2,y_2)$ is bigger than the element at $(z_1,y_1)$. Then the Frobenius number is given by 
$$
g(\x,\y,z)=z_2(s^2+t^2)+y_2 t(3 s^2-t^2)-s(s^2-3 t^2)\,. 
$$

\subsection{The case $s>(2+\sqrt{3})t$}  

If $s>(2+\sqrt{3})t$, then $\z<\y<\x$. The elements of the ($0$-)Ap\'ery set are given as in Table \ref{tb:223-1}, where each point $(Y,X)$ corresponds to the expression $Y \y +X \x $ and the area of the ($0$-)Ap\'ery set is equal to $\z=s^2+t^2$.  

\begin{table}[htbp]
  \centering
\scalebox{0.7}{
\begin{tabular}{ccccccc}
\cline{1-2}\cline{3-4}\cline{5-6}\cline{6-7}
\multicolumn{1}{|c}{$(0,0)$}&$\cdots$&$(t-1,0)$&$(t,0)$&$\cdots$&$\cdots$&\multicolumn{1}{c|}{$(s-1,0)$}\\
\multicolumn{1}{|c}{$\vdots$}&&$\vdots$&$\vdots$&&&\multicolumn{1}{c|}{$\vdots$}\\
\multicolumn{1}{|c}{$(0,s-1)$}&$\cdots$&$(t-1,s-1)$&$(t,s-1)$&$\cdots$&$\cdots$&\multicolumn{1}{c|}{$(s-1,s-1)$}\\
\cline{4-5}\cline{6-7}
\multicolumn{1}{|c}{$(0,s)$}&$\cdots$&\multicolumn{1}{c|}{$(t-1,s)$}&&&&\\\multicolumn{1}{|c}{$\vdots$}&&\multicolumn{1}{c|}{$\vdots$}&&&&\\
\multicolumn{1}{|c}{$\vdots$}&&\multicolumn{1}{c|}{$\vdots$}&&&&\\
\multicolumn{1}{|c}{$(0,s+t-1)$}&$\cdots$&\multicolumn{1}{c|}{$(t-1,s+t-1)$}&&&&\\
\cline{1-2}\cline{3-3}
\end{tabular}
}
  \caption{${\rm Ap}_0(\x ,\y ,\z )$ when $s>(2+\sqrt{3})t$}
  \label{tb:223-1}
\end{table}

By $s\y-t\x=s t\z$, we get $s\y \equiv t\x \pmod{\z }$ and $s\y>t\x $. By using an additional relation $t\y+s\x=(s-t)(s+t)\z$, the sequence $\{\ell\y\pmod{\z }\}_{\ell=0}^{\z -1}$ is given by 
\begin{align}
&(Y,0)\quad(Y=0,1,\dots,s-1),(Y,t),\quad(Y=0,1,\dots,s-1),\notag\\
&(Y,2 t)\quad(Y=0,1,\dots,s-1),\dots,\notag\\
&(Y,\fl{s/t}t)\quad(Y=0,1,\dots,s-1),\notag\\
&\bigl(Y,(\fl{s/t}+1)t\bigr)\quad(Y=0,1,\dots,t-1),\notag\\ 
&\bigl(Y,(\fl{s/t}+1)t-s\bigr)\quad(Y=0,1,\dots,s-1),\notag\\ 
&\bigl(Y,(\fl{s/t}+2)t-s\bigr)\quad(Y=0,1,\dots,s-1),\dots,\notag\\ 
&\bigl(Y,(2\fl{s/t}+1)t-s\bigr)\quad(Y=0,1,\dots,t-1),\notag\\ 
&\bigl(Y,(2\fl{s/t}+1)t-2 s\bigr)\quad(Y=0,1,\dots,s-1),\dots,\notag\\ 
&\bigl(Y,(3\fl{s/t}+1)t-2 s\bigr)\quad(Y=0,1,\dots,t-1),\notag\\ 
&\dots,\notag\\ 
&(Y,s)\quad(Y=0,1,\dots,t-1)\,.
\label{sq:223-1}
\end{align}
Note that $\gcd(s,t)=1$. 
Namely, the next point of $(t-1,s)$ by adding $\y \pmod{\z }$ is $(0,0)$. 
As $\gcd(\y,\z)=1$, the sequence $\{\ell\y\pmod{\z }\}_{\ell=0}^{\z -1}$ matches the sequence $\{\ell\pmod{\z }\}_{\ell=0}^{\z -1}$.  

Since $s>(2+\sqrt{3})t$, we get $(s-1)\y+(s-1)\x>(t-1)\y+(s+t-1)\x$. Hence, by Lemma \ref{lem-mp} (\ref{mp-g}), we have 
$$  
g_0(\x ,\y ,\z ) 
=(s-1)\y +(s-1)\x -\z\,. 
$$ 


\section{The case for $r=4$}  
  
Let $r=4$. The triple of the Diophantine equation $x^2+y^2=z^4$ is parameterized as  
\begin{align*}
(x,y,z)&=\bigl(s^4-6 s^2 t^2+t^4,4 s t(s^2-t^2),s^2+t^2\bigr)\\
&=\bigl((s^2+2 s t-t^2)(s^2-2 s t-t^2),4 s t(s-t)(s+t),s^2+t^2\bigr)\,. 
\end{align*} 
For convenience, we put 
\begin{equation}
\x:=(s^2+2 s t-t^2)(s^2-2 s t-t^2),\quad \y:=4 s t(s-t)(s+t),\quad \z:=s^2+t^2\,.
\label{xyz_r=4}
\end{equation}
Since $\x,\y,\z>0$ and $\gcd(\x,\y,\z)=1$, we see that $s>(1+\sqrt{2})t$, $\gcd(s,t)=1$ and $s\not\equiv t\pmod{2}$.

Then the Frobenius number of this triple is given as follows.   

\begin{theorem}  
If $s^2+t^2<\min\{s^4-6 s^2 t^2+t^4,4 s t(s^2-t^2)\}$, then  
\begin{align*}
&g_0\bigl(s^4-6 s^2 t^2+t^4,4 s t(s^2-t^2),s^2+t^2\bigr)\\
&=\begin{cases}
&(s-1)(s^4-6 s^2 t^2+t^4)+4 s t(s-1)(s^2-t^2)-(s^2+t^2)\\
&\qquad\qquad\qquad\text{if $s>(2+\sqrt{3})t$};\\ 
&(2 s+t-1)(s^4-6 s^2 t^2+t^4)+4 s t(t-1)(s^2-t^2)-(s^2+t^2)\\
&\qquad\qquad\qquad\text{if $2.8350\, t<s<(2+\sqrt{3})t$};\\ 
&(s-1)(s^4-6 s^2 t^2+t^4)+4 s t(s-t-1)(s^2-t^2)-(s^2+t^2)\\
&\qquad\qquad\qquad\text{if $2.5855\, t<s<2.8350\, t$};\\ 
&(5 s+2 t-1)(s^4-6 s^2 t^2+t^4)+4 s t(3 s+t-1)(s^2-t^2)-(s^2+t^2)\\
&\qquad\qquad\qquad\text{if $(1+\sqrt{2})t<s<2.5855\, t$}\,.  
\end{cases}
\end{align*}
Here, $2.8350\dots$ is the larger real root of $x^4-4 x^3+6 x^2-8 x+1=0$. $2.5855\dots$ is the largest real root of $2 x^3-3 x^2-6 x+1=0$. 
\label{th:2-2-4-0}
\end{theorem}

\subsection{The case where $s>\left(1+\sqrt{2}+\sqrt{2(2+\sqrt{2})}\right)t$}  

If $\z$ is not the least, for example, if $\z>\x$, then  
$$
\sqrt{\frac{1+6 s^2-\sqrt{1+16 s^2+32 s^4}}{2}}<t<\sqrt{\frac{1+6 s^2+\sqrt{1+16 s^2+32 s^4}}{2}}\,. 
$$ 
However, any integer pair $(s,t)$, satisfying 
$$
\sqrt{\frac{1+6 s^2-\sqrt{1+16 s^2+32 s^4}}{2}}<t<\frac{s}{1+\sqrt{2}}\,,
$$ 
has not been found for $s<10^7$. Even if it exists, the situation is too complicated. It is similar for the case $\z>\y$. Hence, we assume that $\z<\y<\x$ or $\z<\x<\y$.   
 
Since $\x-\y=s^4-4 s^3 t-6 s^2 t^2+4 s t^3+t^4$ and the real roots of the equation $x^4-4 x^3-6 x^2+4 x+1=0$ are given by  
\begin{align*}
&1-\sqrt{2}-\sqrt{2(2-\sqrt{2})}\,(=-1.4966),\quad 1-\sqrt{2}+\sqrt{2(2-\sqrt{2})}\,(=0.6681)\,,\\
&1+\sqrt{2}-\sqrt{2(2+\sqrt{2})}\,(=-0.1989),\quad 1+\sqrt{2}+\sqrt{2(2+\sqrt{2})}\,(=5.0273)\,,
\end{align*}
we have $\x>\y$ if $s>\left(1+\sqrt{2}+\sqrt{2(2+\sqrt{2})}\right)t$; 
$\x<\y$ if $\bigl(1+\sqrt{2}\bigr)t<s<\left(1+\sqrt{2}+\sqrt{2(2+\sqrt{2})}\right)t$. 

\noindent 
$\bullet$ The case where $s>\left(1+\sqrt{2}+\sqrt{2(2+\sqrt{2})}\right)t$

\noindent 
In this case, $\z<\y<\x$. 
Since $s\y-t\x=t(3 s^2-t^2)\z>0$, we have $s\y\equiv t\x\pmod\z$ and $s\y>t\x$. 
Hence, the arrangement of the elements of the Ap\'ery set is the same as in Table \ref{tb:223-1}, where the position $(Y,X)$ corresponds to the expression $Y\y+X\x$. By using an additional relation $t\y+s\x=s(s^2-3 t^2)\z$, the sequence $\{\ell\y\pmod{\z }\}_{\ell=0}^{\z -1}$ is given as in (\ref{sq:223-1}). 
As $\gcd(\y,\z)=1$, the sequence $\{\ell\y\pmod{\z }\}_{\ell=0}^{\z -1}$ matches the sequence $\{\ell\pmod{\z }\}_{\ell=0}^{\z -1}$. 

Finally, compare the two candidates $(s-1,s-1)$ and $(t-1,s+t-1)$ to take the largest value. Since 
$$
(s-1)\y+(s-1)\x-\bigl((t-1)\y+(s+t-1)\x\bigr)=t\bigl(s^3(3 s-4 t)+2 s^2 t^2+t^3(4 s-t)\bigr)>0\,,
$$ 
the element at $(s-1,s-1)$ is the largest in the Ap\'ery set. Therefore, by Lemma \ref{lem-mp} (\ref{mp-g}), we have 
$$  
g_0(\x ,\y ,\z ) 
=(s-1)\y +(s-1)\x -\z\,. 
$$

\subsection{The case where $(1+\sqrt{2})t<s<\left(1+\sqrt{2}+\sqrt{2(2+\sqrt{2})}\right)t$} 

In this case, $\z<\x<\y$. 

\noindent 
$\bullet$ The case where $(2+\sqrt{3})t<s<\left(1+\sqrt{2}+\sqrt{2(2+\sqrt{2})}\right)t$\\
\noindent 
Since $(s+t)\x-(s-t)\y=(s+t)(s^2-4 s t+t^2)\z>0$, we have $(s+t)\x\equiv(s-t)\y\pmod\z$ and $(s+t)\x>(s-t)\y$. Notice that the real roots of $x^2-4 x+1=0$ are $2-\sqrt{3}$ and $2+\sqrt{3}=3.7320$. 
Hence, the arrangement of the elements of the Ap\'ery set is the same as in Table \ref{tb:223-3}, where the position $(X,Y)$ corresponds to the expression $X\x+Y\y$. 
By using an additional relation $s\x+t\y=s(s^2-3 t^2)\z$, the sequence $\{\ell\y\pmod{\z }\}_{\ell=0}^{\z -1}$ is given as in (\ref{sq:223-3}). 
As $\gcd(\y,\z)=1$, the sequence $\{\ell\y\pmod{\z }\}_{\ell=0}^{\z -1}$ matches the sequence $\{\ell\pmod{\z }\}_{\ell=0}^{\z -1}$. 
Since 
$$
(s-1)\x+(s-1)\y-\bigl((s+t-1)\x+(t-1)\y\bigr)=t\bigl(s^3(3 s-4 t)+2 s^2 t^2+t^3(4 s-t)\bigr)>0\,,
$$ 
the element at $(s-1,s-1)$ is the largest in the Ap\'ery set. Therefore, we have 
$$  
g_0(\x ,\y ,\z ) 
=(s-1)\x +(s-1)\y -\z\,. 
$$ 
\bigskip

\noindent 
$\bullet$ The case where $2.5855\, t<s<(2+\sqrt{3})t$\\
\noindent 
In this case, $(s+t)\x\equiv(s-t)\y\pmod\z$ but $(s+t)\x<(s-t)\y$. Nevertheless, by $(2 s+t)\x-(s-2 t)\y=(2 s^3-3 s^2 t-6 s t^2+t^3)\z>0$ we have $(2 s+t)\x\equiv(s-2 t)\y\pmod\z$ and $(2 s+t)\x>(s-2 t)\y$. Notice that the roots of the equation $2 x^3-3 x^2-6 x+1=0$ are $-1.2413$, $0.1557$ and $2.5855$.  
Since 
$$
(2 s+t-1)\x+(t-1)\y-\bigl((s-1)\x+(s-t-1)\y\bigr)=(s+t)(s^4-4 s^3 t+6 s^2 t^2-8 s t^3+t^4)
$$ 
and the real roots of $x^4-4 x^3+6 x^2-8 x+1=0$ are $0.1380$ and $2.8350$, there are two possibilities. See Table \ref{tb:223-4}. If $2.8350\, t<s<(2+\sqrt{3})t$, the element at $(2 s+t-1,t-1)$ is the largest in the Ap\'ery set. If $2.5855\, t<s<2.8350\, t$, the element at $(s-1,s-t-1)$ is the largest.  
Therefore, if $2.8350\, t<s<(2+\sqrt{3})t$, we have 
$$  
g_0(\x ,\y ,\z ) 
=(2 s+t-1)\x+(t-1)\y -\z\,. 
$$ 
If $2.5855\, t<s<2.8350\, t$ (for example, $(s,t)=(14,5),(17,6),(18,7),\dots$), we have 
$$  
g_0(\x ,\y ,\z ) 
=(s-1)\x+(s-t-1)\y -\z\,. 
$$ 
\bigskip

\noindent 
$\bullet$ The case where $(1+\sqrt{2})t<s<2.5855\, t$\\
\noindent  
In this case, $(2 s+t)\x\equiv(s-2 t)\y\pmod\z$ but $(2 s+t)\x<(s-2 t)\y$.  Nevertheless, by $(5 s+2 t)\x-(2 s-5 t)\y=(5 s^3-6 s^2 t-15 s t^2+2 t^3)\z>0$ we have $(5 s+2 t)\x\equiv(2 s-5 t)\y\pmod\z$ and $(5 s+2 t)\x>(2 s-5 t)\y$. 

Notice that the roots of $5 x^3-6 x^2-15 x+2=0$ are $-1.3142$, $0.1275$ and $2.3867$, and $1+\sqrt{2}>2.3867$

\begin{table}[htbp]
  \centering
\scalebox{0.7}{
\begin{tabular}{ccccccc}
\cline{1-2}\cline{3-4}\cline{5-6}\cline{6-7}
\multicolumn{1}{|c}{$(0,0)$}&$\cdots$&$(3 s+t-1,0)$&$(3 s+t,0)$&$\cdots$&$\cdots$&\multicolumn{1}{c|}{$(5 s+2 t-1,0)$}\\
\multicolumn{1}{|c}{$\vdots$}&&$\vdots$&$\vdots$&&&\multicolumn{1}{c|}{$\vdots$}\\
\multicolumn{1}{|c}{$(0,3 t-s-1)$}&$\cdots$&$(3 s+t-1,3 t-s-1)$&$(3 s+t,3 t-s-1)$&$\cdots$&$\cdots$&\multicolumn{1}{c|}{$(5 s+2 t-1,3 t-s-1)$}\\
\cline{4-5}\cline{6-7}
\multicolumn{1}{|c}{$(0,3 t-s)$}&$\cdots$&\multicolumn{1}{c|}{$(3 s+t-1,3 t-s)$}&&&&\\\multicolumn{1}{|c}{$\vdots$}&&\multicolumn{1}{c|}{$\vdots$}&&&&\\
\multicolumn{1}{|c}{$\vdots$}&&\multicolumn{1}{c|}{$\vdots$}&&&&\\
\multicolumn{1}{|c}{$(0,s-2 t-1)$}&$\cdots$&\multicolumn{1}{c|}{$(3 s+t-1,s-2 t-1)$}&&&&\\
\cline{1-2}\cline{3-3}
\end{tabular}
}
  \caption{${\rm Ap}_0(\x ,\y ,\z )$ when $2.3867\, t<s<2.5855\, t$}
  \label{tb:224-4}
\end{table}

Since 
\begin{align*}
&(5 s+2 t-1)\x+(3 t-s-1)\y-\bigl((3 s+t-1)\x+(s-2 t-1)\y\bigr)\\
&=(s+t)(s^4-4 s^3 t+6 s^2 t^2-8 s t^3+t^4)>0\,,
\end{align*} 
the element at $(5 s+2 t-1,s-1)$ is the largest in the Ap\'ery set. 
Notice that the real roots of $2 x^5-7 x^4+8 x^3+2 x^2-18 x+1=0$ are $-1.0842$, $0.0559$ and $2.2938$. 
Therefore, by Lemma \ref{lem-mp} (\ref{mp-g}), we have 
$$  
g_0(\x ,\y ,\z ) 
=(5 s+2 t-1)\y +(3 s+t-1)\x -\z\,. 
$$ 
\bigskip

\section{The case for $r=5$}  
  
Let $r=5$. The triple of the Diophantine equation $x^2+y^2=z^5$ is parameterized as  
\begin{equation} 
(x,y,z)=\bigl(s(s^4-10 s^2 t^2+5 t^4),t(5 s^4-10 s^2 t^2+t^4),s^2+t^2\bigr)\,. 
\label{xyz_r=5} 
\end{equation}  
For convenience, we put 
$$
\x:=s(s^4-10 s^2 t^2+5 t^4),\quad \y:=t(5 s^4-10 s^2 t^2+t^4),\quad \z:=s^2+t^2\,. 
$$ 
Since $\x,\y,\z>0$ and $\gcd(\x,\y,\z)=1$, we see that $\gcd(s,t)=1$, $s\not\equiv t\pmod{2}$, and 
$$
0<s<\sqrt{\frac{5-2\sqrt{5}}{5}}\,t=0.3249\,t\quad\hbox{or}\quad s>\sqrt{5+2\sqrt{5}}\,t=3.0776\, t\,. 
$$   
Notice that $0.3249$ is the smaller positive root of $5 x^4-10 x^2+1=0$ and $3.0776$ is the larger positive root of $x^4-10 x^2+5=0$.

Then the Frobenius number of this triple is given as follows.   

\begin{theorem}   
\begin{align*}
&g_0\bigl(s(s^4-10 s^2 t^2+5 t^4),t(5 s^4-10 s^2 t^2+t^4),s^2+t^2\bigr)\\
&=\begin{cases}
&(s-1)\x+(s-1)\y-\z \\
&\qquad\qquad\qquad\text{if $s>\left(1+\sqrt{2}+\sqrt{2(2+\sqrt{2})}\right)t$}\\ 
&\qquad\qquad\qquad\text{or $\left(1+\sqrt{5}-\sqrt{5+2\sqrt{5}}\right)t<s<\left(-1-\sqrt{2}+\sqrt{2(2+\sqrt{2})}\right)t$};\\ 
&(2 s+t-1)\x +(t-1)\y -\z\\
&\qquad\qquad\qquad\text{if $4.1894\, t<s<\left(1+\sqrt{2}+\sqrt{2(2+\sqrt{2})}\right)t$};\\
&(s-1)\x+(s-t-1)\y -\z \\
&\qquad\qquad\qquad\text{if $\left(1+\sqrt{5}+\sqrt{2(5+\sqrt{5})}\right)t/2<s<4.1894\, t$};\\ 
&(s-1)\x+(s-2 t-1)-\z \\
&\qquad\qquad\qquad\text{if $\left(1+\sqrt{10}+\sqrt{2(10+\sqrt{10}}\right)t/3<s<3.1098\, t$}\,,\\
&(3 s+t-1)\x+(t-1)\y-\z \\
&\qquad\qquad\qquad\text{if $3.1098\, t<s<\left(1+\sqrt{5}+\sqrt{2(5+\sqrt{5})}\right)t/2$}\,,\\
&(s^2+t^2-1)\x-\z \\
&\qquad\qquad\qquad\text{if $\sqrt{5+2\sqrt{5}}\,t<s<\left(1+\sqrt{10}+\sqrt{2(10+\sqrt{10}}\right)t/3$}\\
&\qquad\qquad\qquad\text{or $\left(-1-\sqrt{5}+\sqrt{2(5+\sqrt{5})}\right)t/2<s<\sqrt{1-2/\sqrt{5}}\, t$}\,,\\
\end{cases}
\end{align*}
\begin{align*}
&g_0\bigl(s(s^4-10 s^2 t^2+5 t^4),t(5 s^4-10 s^2 t^2+t^4),s^2+t^2\bigr)\\
&=\begin{cases}
&(s-1)\x+(s+2 t-1)\y-\z \\
&\qquad\qquad\qquad\text{if $\left(-1-\sqrt{2}+\sqrt{2(2+\sqrt{2})}\right)t<s<0.2386\, t$}\,,\\
&(t-s-1)\x+(t-1)\y-\z \\
&\qquad\qquad\qquad\text{if $0.2386\, t<s<\left(-1-\sqrt{5}+\sqrt{2(5+\sqrt{5})}\right)t/2$}\,,\\
&(s-1)\x+(s+3 t-1)\y-\z \\
&\qquad\qquad\qquad\text{if $\left(-1-\sqrt{5}+\sqrt{2(5+\sqrt{5})}\right)t/2<s<0.3215\, t$}\,,\\
&(t-2 s-1)\x+(t-1)\y-\z \\
&\qquad\qquad\qquad\text{if $0.3215\, t<s<\left(-1+\sqrt{10}+\sqrt{2(10-\sqrt{10})}\right)t/3$}\,,\\
&(t-1)\x+(t-1)\y-\z \\
&\qquad\qquad\qquad\text{if $0<s<\left(1+\sqrt{5}-\sqrt{5+2\sqrt{5}}\right)t$}\,.\\  
\end{cases}
\end{align*}
Here, $4.1894\dots$ is the larger positive root of $x^6-4 x^5-15 x^2+4 x+2=0$. $3.1098$ is the larger positive root of $2 x^6-4 x^5-5 x^4-20 x^2+4 x+3=0$. $0.2386$ is the smaller positive root of $2 x^6+4 x^5-15 x^4-4 x+1=0$.  $0.3215$ is the smaller positive root of $3 x^6+4 x^5-20 x^4-5 x^2-4 x+2=0$. 
\label{th:2-2-5-0}
\end{theorem}

\subsection{The case where $s>\left(1+\sqrt{5}+\sqrt{5+2\sqrt{5}}\right)t$}  

We assume that $\z<\y<\x$ or $\z<\x<\y$. 
Since 
$$
\x-\y=s^5-5 s^4 t-10 s^3 t^2+10 s^2 t^3+5 s t^4- t^5
$$
and the real roots of the equation $x^5-5 x^4-10 x^3+10 x^2+5 x-1=0$ are given by 
\begin{align*}
&1,\quad 1-\sqrt{5}-\sqrt{5-2\sqrt{5}}=-1.9626,\quad 1-\sqrt{5}+\sqrt{5-2\sqrt{5}}=-0.5095,\\
&1+\sqrt{5}-\sqrt{5+2\sqrt{5}}=0.1583,\quad 1+\sqrt{5}+\sqrt{5+2\sqrt{5}}=6.3137\,,
\end{align*} 
we obtain that $\z<\y<\x$ if and only if 
$$
0.1583 t=\left(1+\sqrt{5}-\sqrt{5+2\sqrt{5}}\right)t<s<\sqrt{\frac{5-2\sqrt{5}}{5}}\,t=0.3249\, t
$$
or 
$$
s>\left(1+\sqrt{5}+\sqrt{5+2\sqrt{5}}\right)t=6.3137\, t\,; 
$$
$\z<\x<\y$ if and only if 
$$
0<s<\left(1+\sqrt{5}-\sqrt{5+2\sqrt{5}}\right)t=0.1583\, t
$$
or 
$$  
3.0776\, t=\sqrt{5+2\sqrt{5}}\,t<s<\left(1+\sqrt{5}+\sqrt{5+2\sqrt{5}}\right)t=6.3137\, t\,. 
$$

\noindent 
$\bullet$ The case where $s>\left(1+\sqrt{5}+\sqrt{5+2\sqrt{5}}\right)t$ 
\noindent 

In this case, $\z<\y<\x$. Since $s\y-t\x=s t(s-t)(s+t)\z>0$, we get $s\y\equiv t\x\pmod\z$ and $s\y>t\x$. 

In addition, from $t\y+s\x=(s^2+2 s t-t^2)(s^2-2 s t-t^2)\z>0$, we have $t\y\equiv(s+t-j)\x\pmod\z$ and $t\y>(s+t-j)\x$ ($j=1,2,\dots,t$). 
Hence, the arrangement of the elements of the Ap\'ery set is the same as in Table \ref{tb:223-1}, where the position $(Y,X)$ corresponds to the expression $Y\y+X\x$. Hence, the sequence $\{\ell\y\pmod{\z }\}_{\ell=0}^{\z -1}$ is given as in (\ref{sq:223-1}). 
As $\gcd(\y,\z)=1$, the sequence $\{\ell\y\pmod{\z }\}_{\ell=0}^{\z -1}$ matches the sequence $\{\ell\pmod{\z }\}_{\ell=0}^{\z -1}$. 

Finally, compare the two candidates $(s-1,s-1)$ and $(t-1,s+t-1)$ to take the largest value. Since 
$$
(s-1)\y+(s-1)\x-\bigl((t-1)\y+(s+t-1)\x\bigr)=t(4 s^5-5 s^4 t+10 s^2 t^3-4 s t^4-t^5)>0\,,
$$ 
the element at $(s-1,s-1)$ is the largest in the Ap\'ery set. Notice that the largest real root of the equation $4 x^5-5 x^4+10 x^2-4 x-1=0$ is $0.621267$, which is less than $1+\sqrt{5}+\sqrt{5+2\sqrt{5}}$. Therefore, by Lemma \ref{lem-mp} (\ref{mp-g}), we have 
$$  
g_0(\x ,\y ,\z ) 
=(s-1)\y +(s-1)\x -\z\,. 
$$

\subsection{The case where $\sqrt{5+2\sqrt{5}}\,t<s<\left(1+\sqrt{5}+\sqrt{5+2\sqrt{5}}\right)t$}  

In this case, $\z<\x<\y$. 

\noindent 
$\bullet$ The case where $\left(1+\sqrt{2}+\sqrt{2(2+\sqrt{2})}\right)t<s<\left(1+\sqrt{5}+\sqrt{5+2\sqrt{5}}\right)t$  
\noindent 

For example, $(s,t)=(21,4)$ satisfies this case.    
Since $(s+t)\x-(s-t)\y=(s^4-4 s^3 t-6 s^2 t^2+4 s t^3+t^4)\z$ and the positive real roots of $x^4-4 x^3-6 x^2+4 x+1=0$ are 
$$
1-\sqrt{2}+\sqrt{2(2-\sqrt{2})}=0.6681\quad\hbox{and}\quad 1+\sqrt{2}+\sqrt{2(2+\sqrt{2})}=5.02734\,, 
$$ 
if $\left(1+\sqrt{2}+\sqrt{2(2+\sqrt{2})}\right)t<s<\left(1+\sqrt{5}+\sqrt{5+2\sqrt{5}}\right)t$, then $(s+t)\x\equiv(s-t)\y\pmod\z$ and $(s+t)\x>(s-t)\y$. 
Then the area of the elements of the Ap\'ery set is the same as in Table \ref{tb:223-4}, where the position $(X,Y)$ corresponds to the expression $X\x+Y\y$. 

Compare the elements at $(s-1,s-1)$ and $(s+t-1,t-1)$. We see that  
$$
(s-1)\x+(s-1)\y-\bigl((s+t-1)\x+(t-1)\y\bigr)=t(4 s^5-5 s^4 t+10 s^2 t^3-4 s t^4-t^5)>0
$$ 
because only the positive real root of $4 x^5-5 x^4+10 x^2-4 x-1=0$ is $0.6212$, which is less than $1+\sqrt{2}+\sqrt{2(2+\sqrt{2})}=5.0273$. 
Then the element at $(s-1,s-1)$ is the largest in the Ap\'ery set, and we have 
$$ 
g\bigl(\x ,\y ,\z \bigr)
=(s-1)\x+(s-1)\y -\z\,. 
$$

\noindent 
$\bullet$ The case where $\left(1+\sqrt{5}+\sqrt{2(5+\sqrt{5})}\right)t/2<s<\left(1+\sqrt{2}+\sqrt{2(2+\sqrt{2})}\right)t$  
\noindent 

If $3.5201\, t<s<5.0273\, t$, then $(s+t)\x\equiv(s-t)\y\pmod\z$ but $(s+t)\x<(s-t)\y$. Nevertheless, $(2 s+t)\x-(s-2 t)\y=2(s^4-2 s^3 t-6 s^2 t^2+2 s t^3+t^4)>0$ because the real positive roots of $x^4-2 x^3-6 x^2+2 x+1=0$ are 
$$ 
\frac{1-\sqrt{5}+\sqrt{2(5-\sqrt{5})}}{2}=0.5575\quad\hbox{and}\quad \frac{1+\sqrt{5}+\sqrt{2(5+\sqrt{5})}}{2}=3.5201\,. 
$$ 
Then, $(2 s+t)\x\equiv(s-2 t)\y\pmod\z$ and $(2 s+t)\x>(s-2 t)\y$. Hence, the arrangement of all the elements of the ($0$-)Ap\'ery set are the same as in Table \ref{tb:223-4}. 

Compare the elements at $(2 s+t-1,t-1)$ and $(s-1,s-t-1)$. Since 
$$
(2 s+t-1)\x+(t-1)\y-\bigl((s-1)\x+(s-t-1)\y\bigr)=s^6-4 s^5 t-15 s^2 t^4+4 s t^5+2 t^6
$$ 
and the positive real roots of $x^6-4 x^5-15  x^2+4 x+2=0$ are $0.5114$ and $4.1894$, if 
$$
4.1894\, t<s<\left(1+\sqrt{2}+\sqrt{2(2+\sqrt{2})}\right)t=6.3137\, t 
$$ 
(for example, $(s,t)=(15,4)$), then the element at $(2 s+t-1,t-1)$ is the largest in the Ap\'ery set, and we have 
$$ 
g\bigl(\x ,\y ,\z \bigr)
=(2 s+t-1)\x +(t-1)\y -\z\,. 
$$    
If 
$$
3.0776 t=\sqrt{5+2\sqrt{5}}\,t<s<4.1894 t 
$$ 
(for example, $(s,t)=(17,4)$), then the element at $(2 s+t-1,t-1)$ is the largest, and we have 
$$ 
g\bigl(\x ,\y ,\z \bigr)
=(s-1)\x+(s-t-1)\y -\z\,. 
$$

\noindent 
$\bullet$ The case where $\left(1+\sqrt{10}+\sqrt{2(10+\sqrt{10}}\right)t/3<s$\\$<\left(1+\sqrt{5}+\sqrt{2(5+\sqrt{5})}\right)t/2$  
\noindent 
 
If $3.0976\, t<s<3.5201\, t$, then $(2 s+t)\x\equiv(s-2 t)\y\pmod\z$ but $(2 s+t)\x<(s-2 t)\y$. Nevertheless, $(3 s+t)\x-(s-3 t)\y=(3 s^4-4 s^3 t-18 s^2 t^2+4 s t^3+3 t^4)\z>0$ because the positive real roots of $3 x^4-4 x^3-18 x^2+4 x+3=0$ are 
$$
\frac{1-\sqrt{10}+\sqrt{2(10-\sqrt{10})}}{3}=0.5119
$$ 
and
$$
\frac{1+\sqrt{10}+\sqrt{2(10+\sqrt{10})}}{3}=3.0976\,.
$$ 
Hence, $(3 s+t)\x\equiv(s-3 t)\y\pmod\z$ and $(3 s+t)\x>(s-3 t)\y$.

\begin{table}[htbp]
  \centering
\scalebox{0.9}{
\begin{tabular}{ccccc}
\cline{1-2}\cline{3-4}\cline{5-5}
\multicolumn{1}{|c}{$(0,0)$}&$\cdots$&$\cdots$&$\cdots$&\multicolumn{1}{c|}{}\\
\multicolumn{1}{|c}{}&$\cdots$&$\cdots$&$\cdots$&\multicolumn{1}{c|}{$(3s+t-1,t-1)$}\\
\cline{4-5}
\multicolumn{1}{|c}{}&$\cdots$&\multicolumn{1}{c|}{$\cdots$}&&\\
\multicolumn{1}{|c}{}&$\cdots$&\multicolumn{1}{c|}{$(s-1,s-2 t-1)$}&&\\
\cline{1-2}\cline{3-3}
\end{tabular}
}
  \caption{${\rm Ap}_0(\x ,\y ,\z )$ when $3.0976\, t<s<3.5201\, t$}
  \label{tb:225-24}
\end{table}

Comparing the possible elements (see Table \ref{tb:225-24}), 
\begin{multline*}
(3 s+t-1)\x+(t-1)\y-\bigl((s-1)\x+(s-2 t-1)\y)\\
=2 s^6-4 s^5 t-5 s^4 t^2-20 s^2 t^4+4 s t^5+3 t^6\,. 
\end{multline*} 
The positive real roots of $2 x^6-4 x^5-5 x^4-20 x^2+4 x+3=0$ are $0.4782$ and $3.1098$.  
Hence, if $3.0976\, t<s<3.1098\, t$ (for example, $(s,t)=(31,10),(87,28)$), then $(s-1,s-2 t-1)$ is the largest and 
$$ 
g\bigl(\x ,\y ,\z \bigr)
=(s-1)\x+(s-2 t-1)-\z\,. 
$$  
If $3.1098\, t<s<3.5201\, t$ (for example, $(s,t)=(10,3),(13,4)$), then $(3 s+t-1,t-1)$ is the largest and 
$$ 
g\bigl(\x ,\y ,\z \bigr)
=(3 s+t-1)\x+(t-1)\y-\z\,. 
$$

\noindent 
$\bullet$ The case where $\sqrt{5+2\sqrt{5}}\,t<s<\left(1+\sqrt{10}+\sqrt{2(10+\sqrt{10}}\right)t/3$  
\noindent 

If $3.0776\, t<s<3.0976\, t$, then $(3 s+t)\x\equiv(s-3 t)\y\pmod\z$ but $(3 s+t)\x<(s-3 t)\y$. 
For example, $(34,11),(37,12)$ are in this case. 
The largest element in the Ap\'ery set is $(s^2+t^2-1,0)$. The pattern is similar to that in Table \ref{tb:223-11}. Note that $(4 s+t)\x\equiv(s-4 t)\y\pmod\z$ but $s-4 t<0$.  
Hence,   
$$ 
g\bigl(\x ,\y ,\z \bigr)
=(s^2+t^2-1)\x-\z\,. 
$$

\subsection{The case where $\left(1+\sqrt{5}-\sqrt{5+2\sqrt{5}}\right)t<s<\sqrt{1-2/\sqrt{5}}\, t$} 

In this case, $\x>\y$.   
As $0.1583 t<s<0.3249 t$, we get $(s,t)=(1,4),(2,7),(2,9)$,\\$(3,10),(2,11),(2,13),(3,14),(4,15),(3,16),(5,16),(4,17),(5,18),(4,19),\dots$. 

\noindent 
$\bullet$ The case where $\left(1+\sqrt{5}-\sqrt{5+2\sqrt{5}}\right)t<s<\left(-1-\sqrt{2}+\sqrt{2(2+\sqrt{2})}\right)t$
\noindent 

If $0.1583\, t<s<0.1989\, t$, then $(s,t)=(2,11),(2,13),(3,16),\dots$. 
Since $(s+t)\y-(t-s)\x=(s^4+4 s^3 t-6 s^2 t^2-4 s t^3+t^4)\z$ and the positive real roots of $x^4+4 x^3-6 x^2-4 x+1=0$ are 
$$
-1-\sqrt{2}+\sqrt{2(2+\sqrt{2})}=0.1989\quad\hbox{and}\quad -1+\sqrt{2}+\sqrt{2(2+\sqrt{2})}=1.4966\,, 
$$
we get $(s+t)\y\equiv(t-s)\x\pmod\z$ and $(s+t)\y>(t-s)\x$. Then the elements of the ($0$-)Ap\'ery set are given as in Table \ref{tb:225-3}, where each point $(Y,X)$ corresponds to the expression $Y\y+X\x$. The situation is similar to that in Table \ref{tb:223-3}. 

\begin{table}[htbp]
  \centering
\scalebox{0.9}{
\begin{tabular}{ccccc}
\cline{1-2}\cline{3-4}\cline{5-5}
\multicolumn{1}{|c}{$(0,0)$}&$\cdots$&$\cdots$&$\cdots$&\multicolumn{1}{c|}{}\\
\multicolumn{1}{|c}{}&$\cdots$&$\cdots$&$\cdots$&\multicolumn{1}{c|}{$(s+t-1,s-1)$}\\
\cline{4-5}
\multicolumn{1}{|c}{}&$\cdots$&\multicolumn{1}{c|}{$\cdots$}&&\\
\multicolumn{1}{|c}{}&$\cdots$&\multicolumn{1}{c|}{$(t-1,t-1)$}&&\\
\cline{1-2}\cline{3-3}
\end{tabular}
}
  \caption{${\rm Ap}_0(\x ,\y ,\z )$ when $0.1583\, t<s<0.1989\, t$}
  \label{tb:225-3}
\end{table}

We find that the element at $(t-1,t-1)$ is the largest in the Ap\'ery set because 
$$
(t-1)\y+(t-1)\x-\bigl((s+t-1)\y+(s-1)\x)=-s(s^5+4 s^4 t-10 s^3 t^2+5 s t^4-4 t^5)>0\,. 
$$

Note that only the positive real root of $x^5+4 x^4-10 x^3+5 x-4=0$ is $1.6096$, which is larger than $0.1989$. Hence,   
$$ 
g\bigl(\x ,\y ,\z \bigr)
=(t-1)\y+(t-1)\x-\z\,. 
$$

\noindent 
$\bullet$ The case where $-1-\sqrt{2}+\sqrt{2(2+\sqrt{2})}\,t<s<\left(-1-\sqrt{5}+\sqrt{2(5+\sqrt{5})}\right)t/2$
\noindent 

If $0.1989\, t<s<0.2840\, t$, then $(s,t)=(2,9),(3,14),(4,15),(18,5),\dots$.  
In this case, $(s+t)\y\equiv(t-s)\x\pmod\z$ but $(s+t)\y<(t-s)\x$. Nevertheless, by $(s+2 t)\y-(t-2 s)\x=2(s^4+2 s^3 t-6 s^2 t^2-2 s t^3+t^4)\z$ and the positive real roots of $x^4+2 x^3-6 x^2-2 x+1=0$ are 
$$
\frac{-1-\sqrt{5}+\sqrt{2(5+\sqrt{5})}}{2}=0.2840\quad\hbox{and}\quad \frac{-1+\sqrt{5}+\sqrt{2(5-\sqrt{5})}}{2}=1.7936\,, 
$$
we get 
$(s+2 t)\y\equiv(t-2 s)\x\pmod z$ and $(s+2 t)\y>(t-2 s)\x$.  Comparing the possible largest values of the Ap\'ery set, 
$$
(t-s-1)\y+(t-1)\x-\bigl((s+2 t-1)\y+(s-1)\x)=-(2 s^6+4 s^5 t-15 s^4 t^2-4 s t^5+t^6)
$$
and the positive real roots of $2 x^6+4 x^5-15 x^4-4 x+1=0$ are $0.2386$ and $1.9552$. 

\begin{table}[htbp]
  \centering
\scalebox{0.9}{
\begin{tabular}{ccccc}
\cline{1-2}\cline{3-4}\cline{5-5}
\multicolumn{1}{|c}{$(0,0)$}&$\cdots$&$\cdots$&$\cdots$&\multicolumn{1}{c|}{}\\
\multicolumn{1}{|c}{}&$\cdots$&$\cdots$&$\cdots$&\multicolumn{1}{c|}{$(s+2 t-1,s-1)$}\\
\cline{4-5}
\multicolumn{1}{|c}{}&$\cdots$&\multicolumn{1}{c|}{$\cdots$}&&\\
\multicolumn{1}{|c}{}&$\cdots$&\multicolumn{1}{c|}{$(t-1,t-s-1)$}&&\\
\cline{1-2}\cline{3-3}
\end{tabular}
}
  \caption{${\rm Ap}_0(\x ,\y ,\z )$ when $0.1989\, t<s<0.2840\, t$}
  \label{tb:225-4}
\end{table}
 
See Table \ref{tb:225-4}.  
Hence, if $0.1989\, t<s<0.2386\, t$ (for example, $(s,t)=(1,4),(2,9),(3,14),\dots$), then $(s+2 t-1,s-1)$ takes the largest element and 
$$ 
g\bigl(\x ,\y ,\z \bigr)
=(s+2 t-1)\y+(s-1)\x-\z\,. 
$$   
If $0.2386\, t<s<0.2840\, t$ (for example, $(s,t)=(1,4),(4,15),(18,5),\dots$), then $(t-1,t-s-1)$ takes the largest element and  
$$ 
g\bigl(\x ,\y ,\z \bigr)
=(t-1)\y+(t-s-1)\x-\z\,. 
$$

\noindent 
$\bullet$ The case where $\left(-1-\sqrt{5}+\sqrt{2(5+\sqrt{5})}\right)t/2<s$\\$<\left(-1+\sqrt{10}+\sqrt{2(10-\sqrt{10})}\right)t/3$
\noindent 

If $0.2840 t<s<0.3228 t$, then $(s,t)=(4,13),(5,16),(6,19),(10,31)\dots$. In this case, $(s+2 t)\y\equiv(t-2 s)\x\z$ but $(s+2 t)\y<(t-2 s)\x$. Nevertheless, by $(s+3 t)\y-(t-3 s)\x=(3 s^4+4 s^3 t-18 s^2 t^2-4 s t^3+3 t^4)\z$ and the positive real roots of $3 x^4+4 x^3-18 x^2-4 x+3=0$ are  
$$
\frac{-1-\sqrt{10}+\sqrt{2(10+\sqrt{10})}}{3}=1.9534
$$
and
$$
\frac{-1+\sqrt{10}+\sqrt{2(10-\sqrt{10})}}{3}=0.3228\,, 
$$ 
we get 
$(s+3 t)\y\equiv(t-3 s)\x\pmod\z$ and $(s+3 t)\y>(t-3 s)\x$. Comparing the possible largest values of the Ap\'ery set, 
\begin{multline*}
(t-1)\y+(t-2 s-1)\x-\bigl((s+3 t-1)\y+(s-1)\x)\\
=-(3 s^6+4 s^5 t-20 s^4 t^2-5 s^2 t^4-4 s t^5+2 t^6)
\end{multline*} 
and the positive real roots of $3 x^6+4 x^5-20 x^4-5 x^2-4 x+2=0$ are $0.3215$ and $2.0907$.  

\begin{table}[htbp]
  \centering
\scalebox{0.9}{
\begin{tabular}{ccccc}
\cline{1-2}\cline{3-4}\cline{5-5}
\multicolumn{1}{|c}{$(0,0)$}&$\cdots$&$\cdots$&$\cdots$&\multicolumn{1}{c|}{}\\
\multicolumn{1}{|c}{}&$\cdots$&$\cdots$&$\cdots$&\multicolumn{1}{c|}{$(s+3 t-1,s-1)$}\\
\cline{4-5}
\multicolumn{1}{|c}{}&$\cdots$&\multicolumn{1}{c|}{$\cdots$}&&\\
\multicolumn{1}{|c}{}&$\cdots$&\multicolumn{1}{c|}{$(t-1,t-2 s-1)$}&&\\
\cline{1-2}\cline{3-3}
\end{tabular}
}
  \caption{${\rm Ap}_0(\x ,\y ,\z )$ when $0.2840\, t<s<0.3228\, t$}
  \label{tb:225-5}
\end{table}

See Table \ref{tb:225-5}. 
Hence, if $0.2840\, t<s<0.3215\, t$ (for example, $(s,t)=(4,13),(5,16),(6,19)$), then $(s+3 t-1,s-1)$ takes the largest element and 
$$ 
g\bigl(\x ,\y ,\z \bigr)
=(s+3 t-1)\y+(s-1)\x-\z\,. 
$$   
If $0.3215\, t<s<0.3228\, t$ (for example, $(s,t)=(10,31)$), then $(t-1,t-2 s-1)$ takes the largest element and 
$$ 
g\bigl(\x ,\y ,\z \bigr)
=(t-1)\y+(t-2 s-1)\x-\z\,. 
$$

\noindent 
$\bullet$ The case where $\left(-1+\sqrt{10}+\sqrt{2(10-\sqrt{10})}\right)t/3<s<\sqrt{1-2/\sqrt{5}}\,t$
\noindent 

If $0.3228\, t<s<0.3249\, t$, then $(s,t)=(11,34),(12,37),\dots$. In this case, $(s^2+t^2-1,0)$ takes the largest in the Ap\'ery set (see Table \ref{tb:225-6}) and 
$$ 
g\bigl(\x ,\y ,\z \bigr)
=(s^2+t^2-1)\y-\z\,. 
$$  

\begin{table}[htbp]
  \centering
\begin{tabular}{ccc}
\cline{1-2}\cline{3-3}
\multicolumn{1}{|c}{$(0,0)$}&$\cdots$&\multicolumn{1}{c|}{$(s^2+t^2-1,0)$}\\
\cline{1-2}\cline{3-3}
\end{tabular}
  \caption{${\rm Ap}_0(\x ,\y ,\z )$ when $0.3228\, t<s<0.3249\, t$}
  \label{tb:225-6}
\end{table}

\subsection{The case where $0<s<\left(1+\sqrt{5}-\sqrt{5+2\sqrt{5}}\right)t$} 

In this case, $\x<\y$.   
Since $0<s<0.1538\, t\, t$, we can get $(s,t)=(1,8+2 j)$ ($j=0,1,\dots$), $(2,15+2 j)$ ($j=0,1,\dots$), $(3,20),(3,22),(3,26),\dots$. 
Since $t\x-s\y=4 s t(t-s)(s+t)\z$, we have $t\x\equiv s\y\pmod\z$ and $t\x>s\y$.  

Then the elements of the ($0$-)Ap\'ery set are given as in Table \ref{tb:225-12}, where each point $(X,Y)$ corresponds to the expression $X\x+Y\y$. The situation is similar to that in Table \ref{tb:225-12}.  

\begin{table}[htbp]
  \centering
\scalebox{0.9}{
\begin{tabular}{ccccc}
\cline{1-2}\cline{3-4}\cline{5-5}
\multicolumn{1}{|c}{$(0,0)$}&$\cdots$&$\cdots$&$\cdots$&\multicolumn{1}{c|}{}\\
\multicolumn{1}{|c}{}&$\cdots$&$\cdots$&$\cdots$&\multicolumn{1}{c|}{$(t-1,t-1)$}\\
\cline{4-5}
\multicolumn{1}{|c}{}&$\cdots$&\multicolumn{1}{c|}{$\cdots$}&&\\
\multicolumn{1}{|c}{}&$\cdots$&\multicolumn{1}{c|}{$(s-1,s+t-1)$}&&\\
\cline{1-2}\cline{3-3}
\end{tabular}
}
  \caption{${\rm Ap}_0(\x ,\y ,\z )$ when $0<s<\left(1+\sqrt{5}-\sqrt{5+2\sqrt{5}}\right)t$}
  \label{tb:225-12}
\end{table}

Comparing the possible elements, we have 
$$
(t-1)\x+(t-1)\y-\bigl((s-1)\x+(s+t-1)\y\bigr)=-s(s^5+4 s^4 t-10 s^3 t^2+5 s t^4-4 t^5)>0\,.
$$ 
Notice that only the positive real root of $x^5+4 x^4-10 x^3+5 x-4=0$ is $1.6096$. 
Then $(t-1,t-1)$ takes the largest element and 
$$ 
g\bigl(\x ,\y ,\z \bigr)
=(t-1)\x+(t-1)\y-\z\,. 
$$

\section{Sylvester number (genus)}  

\subsection{The case where $r=3$}  

Remember that for convenience, we use $(\x,\y,\z)$ as in (\ref{xyz_r=3}).  

\begin{theorem}
\begin{align*}
&n_0\bigl(s(s^2-3 t^2),t(3 s^2-t^2),s^2+t^2\bigr)\\
&=\begin{cases}
&\frac{1}{2}\bigl((s-1)\x+(t-1)\y+3 s t(s^2-s t-t^2)+(s^2-1)(t^2-1)\bigr)\\
&\qquad\qquad\qquad\text{if $s>(1+\sqrt{2})t$};\\ 
&\frac{1}{2}\bigl((s-1)\x+(t-1)\y+s t(5 s^2-7 s t-5 t^2)+(s^2-1)(t^2-1)\bigr)\\
&\qquad\qquad\qquad\text{if $2 t<s<(1+\sqrt{2})t$};\\ 
&\frac{1}{2}\bigl((s+1)\x-(5 t-1)\y-s t(11 s^2-21 s t-15 t^2)+(s^2-1)(t^2-1)\bigr))\\
&\qquad\qquad\qquad\text{if $(2+\sqrt{13})t/3<s<2 t$};\\ 
&\frac{1}{2}\bigl((13 s+1)\x-(17 t-1)\y-s t(45 s^2-93 s t-49 t^2)+(s^2-1)(t^2-1)\bigr)\\
&\qquad\qquad\qquad\text{if $(3+\sqrt{34})t/5<s<(2+\sqrt{13})t/3$};\\ 
&\frac{1}{2}\bigl((33 s+1)\x-(37 t-1)\y-s t(99 s^2-209 s t-103 t^2)+(s^2-1)(t^2-1)\bigr)\\
&\qquad\qquad\qquad\text{if $\sqrt{3}\,t<s<(3+\sqrt{34})t/5$}\,.  
\end{cases}
\end{align*}
\label{th:2-2-3-1}
\end{theorem}   

\noindent 
$\bullet$ The case where $(1+\sqrt{2})\,t<s<(2+\sqrt{3})\,t$  

By using Table \ref{tb:223-3}, we have 
\begin{align*}  
\sum_{w\in{\rm Ap}(\x,\y,\z)}w&=\sum_{X=0}^{s+t-1}\sum_{Y=0}^{t-1}(X\x+Y\y)+\sum_{X=0}^{s-1}\sum_{Y=t}^{s-1}(X\x+Y\y)\\
&=\frac{\z}{2}(-s^3+s^4-3 s^2 t+3 s^3 t+3 s t^2-2 s^2 t^2+t^3-3 s t^3-t^4)\,. 
\end{align*}

By Lemma \ref{lem-mp} (\ref{mp-n}), 
\begin{align*}  
n_0(\x,\y,\z)&=\frac{1}{\z}\sum_{w\in{\rm Ap}(\x,\y,\z)}w-\frac{\z-1}{2}\\
&=\frac{1}{2}\bigl((s-1)\x+(t-1)\y+3 s t(s^2-s t-t^2)+(s^2-1)(t^2-1)\bigr)\,. 
\end{align*}

\noindent 
$\bullet$ The case where $2 t<s<(1+\sqrt{2})\,t$  

\begin{align*}  
&n_0(\x,\y,\z)\\
&=\frac{1}{2}\bigl((s-1)\x+(t-1)\y+s t(5 s^2-7 s t-5 t^2)+(s^2-1)(t^2-1)\bigr)\,.
\end{align*}

\noindent 
$\bullet$ The case where $(2+\sqrt{13})t/3<s<2 t$  

\begin{align*}  
&n_0(\x,\y,\z)\\
&=\frac{1}{2}\bigl((s+1)\x-(5 t-1)\y-s t(11 s^2-21 s t-15 t^2)+(s^2-1)(t^2-1)\bigr)\,. 
\end{align*}

\noindent 
$\bullet$ The case where $(3+\sqrt{34})t/5<s<(2+\sqrt{13})t/3$  

\begin{align*}  
&n_0(\x,\y,\z)\\
&=\frac{1}{2}\bigl((13 s+1)\x-(17 t-1)\y-s t(45 s^2-93 s t-49 t^2)+(s^2-1)(t^2-1)\bigr)\,. 
\end{align*}

\noindent 
$\bullet$ The case where $\sqrt{3}\,t<s<(3+\sqrt{34})t/5$  

\begin{align*}  
&n_0(\x,\y,\z)\\
&=\frac{1}{2}\bigl((33 s+1)\x-(37 t-1)\y-s t(99 s^2-209 s t-103 t^2)+(s^2-1)(t^2-1)\bigr)\,. 
\end{align*}

\subsubsection{The case where $\x<\z<\y$} 

If the Ap\'ery set is arranged as in Table \ref{tb:223-11}, then 
by 
$$
\sum_{w\in{\rm Ap}_0(A)}w=\sum_{Z=0}^{z_0}Z\z=\frac{z_0(z_0+1)}{2}\z\,,
$$ 
where $z_0+1=\x$, 
we have 
\begin{align*}  
n_0(\x,\y,\z)&=\frac{1}{\x}\sum_{w\in{\rm Ap}(\x,\y,\z)}w-\frac{\x-1}{2}\\
&=\frac{(\x-1)(\z-1)}{2}\,. 
\end{align*}   
If the Ap\'ery set is arranged as in Table \ref{tb:223-12}, then 
by 
$$
\sum_{w\in{\rm Ap}_0(A)}w=\sum_{Z=0}^{z_1}\sum_{Y=0}^{y_1}(Z\z+Y\y)+\sum_{Z=0}^{z_2}\sum_{Y=y_1+1}^{y_2}(Z\z+Y\y)\,, 
$$
we can calculate $n_0(\x,\y,\z)$ individually.  


\subsubsection{The case $s>(2+\sqrt{3})t$}  

By using Table \ref{tb:223-1}, we have 
\begin{align*}  
\sum_{w\in{\rm Ap}(\x,\y,\z)}w&=\sum_{Y=0}^{s-1}\sum_{X=0}^{s-1}(X\x+Y\y)+\sum_{Y=0}^{t-1}\sum_{X=s}^{s+t-1}(X\x+Y\y)\\
&=\frac{\z}{2}(-s^3+s^4-3 s^2 t+3 s^3 t+3 s t^2-2 s^2 t^2+t^3-3 s t^3-t^4)\,. 
\end{align*}
By Lemma \ref{lem-mp} (\ref{mp-n}), 
\begin{align*}  
n_0(\x,\y,\z)&=\frac{1}{\z}\sum_{w\in{\rm Ap}(\x,\y,\z)}w-\frac{\z-1}{2}\\
&=\frac{1}{2}\bigl((s-1)\x+(t-1)\y+3 s t(s^2-s t-t^2)+(s^2-1)(t^2-1)\bigr)\,. 
\end{align*}

\subsection{The case for $r=4$}  
 
Remember that for convenience, we use $(\x,\y,\z)$ as in (\ref{xyz_r=4}).     

\begin{theorem}  
If $s^2+t^2<\min\{s^4-6 s^2 t^2+t^4,4 s t(s^2-t^2)\}$, then  
\begin{align*}
&n_0\bigl(s^4-6 s^2 t^2+t^4,4 s t(s^2-t^2),s^2+t^2\bigr)\\
&=\begin{cases}
&\frac{1}{2}\bigl((s+t-1)\x-\y+s t(3(s^3-t^3)+s t(s-t-1))\\
&\quad +(s^2-1)(t^2-1)\bigr)\qquad\text{if $s>(2+\sqrt{3})t$};\\ 
&\frac{1}{2}\bigl((s+t-1)\x-\y+s t(5 s^3-t^3-s t(5 s+7 t+1))\\
&\quad +(s^2-1)(t^2-1)\bigr)\qquad\text{if $2.5855\, t<s<(2+\sqrt{3})t$};\\ 
&\frac{1}{2}\bigl((7 s-7 t+1)\x+\y+s t(-35 s^3+17 t^3+s t(53 s+65 t-1))\\
&\quad +(s^2-1)(t^2-1)\bigr)\qquad\text{if $(1+\sqrt{2})t<s<2.5855\, t$}\,.  
\end{cases}
\end{align*}
Here, $2.8350\dots$ is the larger real root of $x^4-4 x^3+6 x^2-8 x+1=0$. 
\label{th:2-2-4-1}
\end{theorem}  

\subsubsection{The case where $s>\left(1+\sqrt{2}+\sqrt{2(2+\sqrt{2})}\right)t$}  

By using Table \ref{tb:223-1}, we have 
\begin{align*}  
&\sum_{w\in{\rm Ap}(\x,\y,\z)}w=\sum_{Y=0}^{s-1}\sum_{X=0}^{s-1}(Y\y+X\x)+\sum_{Y=0}^{t-1}\sum_{X=s}^{s+t-1}(Y\y+X\x)\\
&=\frac{\z}{2}(-s^4+s^5-4 s^3 t+4 s^4 t+6 s^2 t^2-5 s^3 t^2+4 s t^3-7 s^2 t^3-t^4-2 s t^4+t^5)\,. 
\end{align*}
By Lemma \ref{lem-mp} (\ref{mp-n}), 
\begin{align*}  
&n_0(\x,\y,\z)=\frac{1}{\z}\sum_{w\in{\rm Ap}(\x,\y,\z)}w-\frac{\z-1}{2}\\
&=\frac{1}{2}\bigl((s+t-1)\x-\y+s t(3(s^3-t^3)+s t(s-t-1))+(s^2-1)(t^2-1)\bigr)\,. 
\end{align*}

\subsubsection{The case where $(1+\sqrt{2})t<s<\left(1+\sqrt{2}+\sqrt{2(2+\sqrt{2})}\right)t$} 

\noindent 
$\bullet$ The case where $(2+\sqrt{3})t<s<\left(1+\sqrt{2}+\sqrt{2(2+\sqrt{2})}\right)t$\\
\noindent 

By using Table \ref{tb:223-3}, we have 
\begin{align*}  
&\sum_{w\in{\rm Ap}(\x,\y,\z)}w=\sum_{X=0}^{s+t-1}\sum_{Y=0}^{t-1}(X\x+Y\y)+\sum_{X=0}^{s-1}\sum_{Y=t}^{s-1}(X\x+Y\y)\\
&=\frac{\z}{2}((-s^4+s^5-4 s^3 t+4 s^4 t+6 s^2 t^2-5 s^3 t^2+4 s t^3\\
&\qquad -7 s^2 t^3-t^4-2 s t^4+t^5))\,. 
\end{align*}
Hence, 
\begin{align*}   
&n_0(\x,\y,\z)=\frac{1}{\z}\sum_{w\in{\rm Ap}(\x,\y,\z)}w-\frac{\z-1}{2}\\
&=\frac{1}{2}\bigl((s+t-1)\x-\y+s t(3(s^3-t^3)+s t(s-t-1))+(s^2-1)(t^2-1)\bigr)\,.
\end{align*}

\noindent 
$\bullet$ The case where $2.5855\, t<s<(2+\sqrt{3})t$\\
\noindent 

By Table \ref{tb:223-4}, we have 
\begin{align*}  
&n_0(\x,\y,\z)\\
&=\frac{1}{2}\bigl((s+t-1)\x-\y+s t(5 s^3-t^3-s t(5 s+7 t+1))+(s^2-1)(t^2-1)\bigr)\,. 
\end{align*}

\noindent 
$\bullet$ The case where $(1+\sqrt{2})t<s<2.5855\, t$\\
\noindent  

By Table \ref{tb:224-4}, we have 
\begin{align*}  
&n_0(\x,\y,\z)\\
&=\frac{1}{2}\bigl((7 s-7 t+1)\x+\y+s t(-35 s^3+17 t^3+s t(53 s+65 t-1))\\
&\qquad\qquad +(s^2-1)(t^2-1)\bigr)\,. 
\end{align*}

\subsection{The case for $r=5$}  
  
Remember that for convenience, we use $(\x,\y,\z)$ as in (\ref{xyz_r=5}).  

\begin{theorem}  
\begin{align*}
&n_0\bigl(s(s^4-10 s^2 t^2+5 t^4),t(5 s^4-10 s^2 t^2+t^4),s^2+t^2\bigr)\\
&=\begin{cases}
&\frac{1}{2}\bigl((s-1)\x+(t-1)\y+s t(5(s^4+t^4)-s t(4 s^2+14 s t-4 t^2+1))\\
&\quad +(s^2-1)(t^2-1)\bigr)\qquad\text{if $s>\left(1+\sqrt{2}+\sqrt{2(2+\sqrt{2})}\right)t$}\,;\\ 
&\frac{1}{2}\bigl((s-1)\x+(t-1)\y+s t(7(s^4+t^4)-s t(12 s^2+26 s t-12 t^2+1))\\
&\quad +(s^2-1)(t^2-1)\bigr)\\
&\qquad\text{if $\left(1+\sqrt{5}+\sqrt{2(5+\sqrt{5})}\right)t/2<s<\left(1+\sqrt{2}+\sqrt{2(2+\sqrt{2})}\right)t$}\,;\\
&\frac{1}{2}\bigl((s-1)\x+(t-1)\y+s t(11(s^4+t^4)-s t(20 s^2+50 s t-20 t^2+1))\\
&\quad +(s^2-1)(t^2-1)\bigr)\\
&\qquad\text{if $\left(1+\sqrt{10}+\sqrt{2(10+\sqrt{10}}\right)t/3<s<\left(1+\sqrt{5}+\sqrt{2(5+\sqrt{5})}\right)t/2$};\\ 
&\frac{(\x-1)(\z-1)}{2}\\
&\qquad\text{if $\sqrt{5+2\sqrt{5}}\,t<s<\left(1+\sqrt{10}+\sqrt{2(10+\sqrt{10}}\right)t/3$}\,;\\
&\frac{1}{2}\bigl((s-1)\x+(t-1)\y+s t(5(s^4+t^4)+s t(4 s^2-14 s t-4 t^2-1))\\
&\quad +(s^2-1)(t^2-1)\bigr)\\
&\qquad\text{if $\left(1+\sqrt{5}-\sqrt{5+2\sqrt{5}}\right)t<s<\left(-1-\sqrt{2}+\sqrt{2(2+\sqrt{2})}\right)t$}\,;\\
\end{cases}
\end{align*}
\begin{align*}
&n_0\bigl(s(s^4-10 s^2 t^2+5 t^4),t(5 s^4-10 s^2 t^2+t^4),s^2+t^2\bigr)\\
&=\begin{cases}
&\frac{1}{2}\bigl((s-1)\x+(t-1)\y+s t(7(s^4+t^4)+s t(12 s^2-26 s t-12 t^2-1))\\
&\quad +(s^2-1)(t^2-1)\bigr) \\
&\qquad\text{if $-1-\sqrt{2}+\sqrt{2(2+\sqrt{2})}\,t<s<\left(-1-\sqrt{5}+\sqrt{2(5+\sqrt{5})}\right)t/2$}\,;\\
&\frac{1}{2}\bigl((s-1)\x+(t-1)\y+s t(11(s^4+t^4)+s t(20 s^2-50 s t-120 t^2-1))\\
&\quad +(s^2-1)(t^2-1)\bigr) \\
&\qquad\text{if $\left(-1-\sqrt{5}+\sqrt{2(5+\sqrt{5})}\right)t/2<s<\left(-1+\sqrt{10}+\sqrt{2(10-\sqrt{10})}\right)t/3$}\,;\\
&\frac{(\y-1)(\z-1)}{2} \\
&\qquad\text{if $\left(-1+\sqrt{10}+\sqrt{2(10-\sqrt{10})}\right)t/3<s<\sqrt{1-2/\sqrt{5}}\,t$}\,;\\
&\frac{1}{2}\bigl((s-1)\x+(t-1)\y+s t(5(s^4+t^4)+s t(4 s^2-14 s t-4 t^2+1))\\
&\quad +(s^2-1)(t^2-1)\bigr)
\qquad\text{if $0<s<\left(1+\sqrt{5}-\sqrt{5+2\sqrt{5}}\right)t$}\,.\\  
\end{cases}
\end{align*} 
\label{th:2-2-5-1}
\end{theorem}  

\subsubsection{The case where $s>\left(1+\sqrt{5}+\sqrt{5+2\sqrt{5}}\right)t$}  

By using Table \ref{tb:223-1}, we have 
\begin{align*}  
\sum_{w\in{\rm Ap}(\x,\y,\z)}w&=\sum_{Y=0}^{s-1}\sum_{X=0}^{s-1}(Y\y+X\x)+\sum_{Y=0}^{t-1}\sum_{X=s}^{s+t-1}(Y\y+X\x)\\
&=\frac{\z}{2}(-s^5+s^6-5 s^4 t+5 s^5 t+10 s^3 t^2-9 s^4 t^2+10 s^2 t^3\\
&\qquad -14 s^3 t^3-5 s t^4-s^2 t^4-t^5+5 s t^5+t^6)\,. 
\end{align*}
By Lemma \ref{lem-mp} (\ref{mp-n}), 
\begin{align*}  
&n_0(\x,\y,\z)=\frac{1}{\z}\sum_{w\in{\rm Ap}(\x,\y,\z)}w-\frac{\z-1}{2}\\
&=\frac{1}{2}\bigl((s-1)\x+(t-1)\y+s t(5(s^4+t^4)-s t(4 s^2+14 s t-4 t^2+1))\\
&\qquad +(s^2-1)(t^2-1)\bigr)\,. 
\end{align*}

\subsubsection{The case where $\sqrt{5+2\sqrt{5}}\,t<s<\left(1+\sqrt{5}+\sqrt{5+2\sqrt{5}}\right)t$}  
 
\noindent 
$\bullet$ The case where $\left(1+\sqrt{2}+\sqrt{2(2+\sqrt{5})}\right)t<s<\left(1+\sqrt{5}+\sqrt{5+2\sqrt{5}}\right)t$  
\noindent 

By using Table \ref{tb:223-3}, we have 
\begin{align*}  
\sum_{w\in{\rm Ap}(\x,\y,\z)}w&=\sum_{X=0}^{s+t-1}\sum_{Y=0}^{t-1}(X\x+Y\y)+\sum_{X=0}^{s-1}\sum_{Y=t}^{s-1}(X\x+Y\y)\\
&=\frac{\z}{2}(-s^5+s^6-5 s^4 t+5 s^5 t+10 s^3 t^2-9 s^4 t^2+10 s^2 t^3\\
&\qquad -14 s^3 t^3-5 s t^4-s^2 t^4-t^5+5 s t^5+t^6)\,. 
\end{align*}
By Lemma \ref{lem-mp} (\ref{mp-n}), 
\begin{align*}  
&n_0(\x,\y,\z)=\frac{1}{\z}\sum_{w\in{\rm Ap}(\x,\y,\z)}w-\frac{\z-1}{2}\\
&=\frac{1}{2}\bigl((s-1)\x+(t-1)\y+s t(5(s^4+t^4)-s t(4 s^2+14 s t-4 t^2+1))\\
&\qquad +(s^2-1)(t^2-1)\bigr)\,. 
\end{align*}

\noindent 
$\bullet$ The case where $\left(1+\sqrt{5}+\sqrt{2(5+\sqrt{5})}\right)t/2<s<\left(1+\sqrt{2}+\sqrt{2(2+\sqrt{2})}\right)t$  
\noindent 

By using Table \ref{tb:223-4}, 
\begin{align*}  
&n_0(\x,\y,\z)\\
&=\frac{1}{2}\bigl((s-1)\x+(t-1)\y+s t(7(s^4+t^4)-s t(12 s^2+26 s t-12 t^2+1))\\
&\qquad +(s^2-1)(t^2-1)\bigr)\,. 
\end{align*}

\noindent 
$\bullet$ The case where $\left(1+\sqrt{10}+\sqrt{2(10+\sqrt{10}}\right)t/3<s$\\$<\left(1+\sqrt{5}+\sqrt{2(5+\sqrt{5})}\right)t/2$  
\noindent 

By using Table \ref{tb:225-24}, 
\begin{align*}  
&n_0(\x,\y,\z)\\
&=\frac{1}{2}\bigl((s-1)\x+(t-1)\y+s t(11(s^4+t^4)-s t(20 s^2+50 s t-20 t^2+1))\\
&\qquad +(s^2-1)(t^2-1)\bigr)\,. 
\end{align*}

\noindent 
$\bullet$ The case where $\sqrt{5+2\sqrt{5}}\,t<s<\left(1+\sqrt{10}+\sqrt{2(10+\sqrt{10}}\right)t/3$  
\noindent 

By using Table \ref{tb:223-11}, 
$$  
n_0(\x,\y,\z)=\frac{(\x-1)(\z-1)}{2}\,. 
$$

\subsubsection{The case where $\left(1+\sqrt{5}-\sqrt{5+2\sqrt{5}}\right)t<s<\sqrt{1-2/\sqrt{5}}\, t$}

\noindent 
$\bullet$ The case where $\left(1+\sqrt{5}-\sqrt{5+2\sqrt{5}}\right)t<s<\left(-1-\sqrt{2}+\sqrt{2(2+\sqrt{2})}\right)t$
\noindent 

By using Table \ref{tb:225-3}, we have 
\begin{align*}  
\sum_{w\in{\rm Ap}(\x,\y,\z)}w&=\sum_{Y=0}^{s-1}\sum_{X=0}^{s-1}(Y\y+X\x)+\sum_{Y=0}^{t-1}\sum_{X=s}^{s+t-1}(Y\y+X\x)\\
&=\frac{\z}{2}(-s^5+s^6-5 s^4 t+5 s^5 t+10 s^3 t^2-s^4 t^2+10 s^2 t^3\\
&\qquad -14 s^3 t^3-5 s t^4-9 s^2 t^4-t^5+5 s t^5+t^6)\,. 
\end{align*}
By Lemma \ref{lem-mp} (\ref{mp-n}), 
\begin{align*}  
&n_0(\x,\y,\z)=\frac{1}{\z}\sum_{w\in{\rm Ap}(\x,\y,\z)}w-\frac{\z-1}{2}\\
&=\frac{1}{2}\bigl((s-1)\x+(t-1)\y+s t(5(s^4+t^4)+s t(4 s^2-14 s t-4 t^2-1))\\
&\qquad +(s^2-1)(t^2-1)\bigr)\,. 
\end{align*}

\noindent 
$\bullet$ The case where $-1-\sqrt{2}+\sqrt{2(2+\sqrt{2})}\,t<s<\left(-1-\sqrt{5}+\sqrt{2(5+\sqrt{5})}\right)t/2$
\noindent 

By using Table \ref{tb:225-4}, we have 
\begin{align*}  
&n_0(\x,\y,\z)\\
&=\frac{1}{2}\bigl((s-1)\x+(t-1)\y+s t(7(s^4+t^4)+s t(12 s^2-26 s t-12 t^2-1))\\
&\qquad +(s^2-1)(t^2-1)\bigr)\,. 
\end{align*}

\noindent 
$\bullet$ The case where $\left(-1-\sqrt{5}+\sqrt{2(5+\sqrt{5})}\right)t/2<s$\\$<\left(-1+\sqrt{10}+\sqrt{2(10-\sqrt{10})}\right)t/3$
\noindent 

By using Table \ref{tb:225-5}, we have 
\begin{align*}  
&n_0(\x,\y,\z)\\
&=\frac{1}{2}\bigl((s-1)\x+(t-1)\y+s t(11(s^4+t^4)+s t(20 s^2-50 s t-120 t^2-1))\\
&\qquad +(s^2-1)(t^2-1)\bigr)\,. 
\end{align*}

\noindent 
$\bullet$ The case where $\left(-1+\sqrt{10}+\sqrt{2(10-\sqrt{10})}\right)t/3<s<\sqrt{1-2/\sqrt{5}}\,t$ 

By using Table \ref{tb:225-6}, we have 
$$
n_0(\x,\y,\z)=\frac{(\y-1)(\z-1)}{2}\,.
$$

\subsubsection{The case where $0<s<\left(1+\sqrt{5}-\sqrt{5+2\sqrt{5}}\right)t$} 

By using Table \ref{tb:225-12}, we have 
\begin{align*}  
\sum_{w\in{\rm Ap}(\x,\y,\z)}w&=\sum_{X=0}^{s+t-1}\sum_{Y=0}^{t-1}(X\x+Y\y)+\sum_{X=0}^{s-1}\sum_{Y=t}^{s-1}(X\x+Y\y)\\
&=\frac{\z}{2}(-s^5+s^6-5 s^4 t+5 s^5 t+10 s^3 t^2-s^4 t^2+10 s^2 t^3\\
&\qquad -14 s^3 t^3-5 s t^4-9 s^2 t^4-t^5+5 s t^5+t^6)\,. 
\end{align*}
By Lemma \ref{lem-mp} (\ref{mp-n}), 
\begin{align*}  
&n_0(\x,\y,\z)=\frac{1}{\z}\sum_{w\in{\rm Ap}(\x,\y,\z)}w-\frac{\z-1}{2}\\
&=\frac{1}{2}\bigl((s-1)\x+(t-1)\y+s t(5(s^4+t^4)+s t(4 s^2-14 s t-4 t^2+1))\\
&\qquad +(s^2-1)(t^2-1)\bigr)\,. 
\end{align*}

\section{$p$-numerical semigroup}  

When $p>0$, the formulas of $p$-Frobenius numbers and $p$-Sylvester number may be obtained though there are many different situations.   
For example, consider the case where $(1+\sqrt{2})\,t<s<(2+\sqrt{3})\,t$ on the Diophantine equation for $r=3$. 
As $(s+t)\x\equiv(s-t)\y\pmod\z$, with respect to the area of the elements of ${\rm Ap}_0(A)$ shown in Table \ref{tb:223-3} (taken as Range 1 (R1) in the longitudinal direction with the side-length $s+t$), the area of the elements of the first ($s-t$) row moves in the form of filling the lower gap (modulo equal) as $0_1\Longrightarrow 1_1$. The portion protruding from the Range 1 is filled further downward from the left in order as $0_2\Longrightarrow 1_2$. The remaining regions except the first ($s-t$) row of the elements of ${\rm Ap}_0(A)$ are shifted (modulo equal) to the upper right side of Range 1 as Range 2 (R2) as $0_3\Longrightarrow 1_3$. It is as shown in Table \ref{tb:223-1-1}. 
Note that $s-t>t$.  

Here, each point $(X,Y)$ corresponds to the expression $X\x+Y\y$ and the area of the $1$-Ap\'ery set is also equal to $\z=s^2+t^2$.   

\begin{table}[htbp]
  \centering
\scalebox{0.7}{
\begin{tabular}{cccccccccccccc}
\multicolumn{1}{|c}{}&&&R1&&&&\multicolumn{1}{|c}{}&&&R2&&&\multicolumn{1}{c|}{}\\
\cline{1-2}\cline{3-4}\cline{5-6}\cline{7-8}\cline{9-10}\cline{11-12}\cline{13-14}
\multicolumn{1}{|c}{}&&\multicolumn{1}{|c}{}&&&&&\multicolumn{1}{|c}{}&&$1_3$&&&\multicolumn{1}{|c}{}&\multicolumn{1}{c|}{}\\
\multicolumn{1}{|c}{}&$0_1$&\multicolumn{1}{|c}{}&&&&&\multicolumn{1}{|c}{}&&&&\ctext{$D$}&\multicolumn{1}{|c}{}&\multicolumn{1}{c|}{}\\
\cline{6-6}\cline{7-8}\cline{9-10}\cline{11-12}
\multicolumn{1}{|c}{}&&\multicolumn{1}{|c}{}&$0_2$&&\multicolumn{1}{|c}{}&&\multicolumn{1}{|c}{}&&&&&&\multicolumn{1}{c|}{}\\
\multicolumn{1}{|c}{}&&\multicolumn{1}{|c}{}&&&\multicolumn{1}{|c}{}&$1_1$&\multicolumn{1}{|c}{}&&&&&&\multicolumn{1}{c|}{}\\
\cline{1-2}\cline{3-4}\cline{5-5}
\multicolumn{1}{|c}{}&&$0_3$&&&\multicolumn{1}{|c}{}&&\multicolumn{1}{|c}{}&&&&&&\multicolumn{1}{c|}{}\\
\multicolumn{1}{|c}{}&&&&&\multicolumn{1}{|c}{}&\ctext{$C$}&\multicolumn{1}{|c}{}&&&&&&\multicolumn{1}{c|}{}\\
\cline{1-2}\cline{3-4}\cline{5-6}\cline{7-7}
\multicolumn{1}{|c}{}&&&&&\multicolumn{1}{|c}{}&&\multicolumn{1}{|c}{}&&&&&&\multicolumn{1}{c|}{}\\
\multicolumn{1}{|c}{}&&&&\ctext{$B$}&\multicolumn{1}{|c}{}&&\multicolumn{1}{|c}{}&&&&&&\multicolumn{1}{c|}{}\\
\cline{4-5}
\multicolumn{1}{|c}{}&&$1_2$&\multicolumn{1}{|c}{}&&&&\multicolumn{1}{|c}{}&&&&&&\multicolumn{1}{c|}{}\\
\multicolumn{1}{|c}{}&&\phantom{ap1-2}\ctext{$A$}&\multicolumn{1}{|c}{}&&&&\multicolumn{1}{|c}{}&&&&&&\multicolumn{1}{c|}{}\\
\cline{1-2}\cline{3-3}
\multicolumn{1}{|c}{}&&&&&&&\multicolumn{1}{|c}{}&&&&&&\multicolumn{1}{c|}{}\\
\multicolumn{1}{|c}{}&&&&&&&\multicolumn{1}{|c}{}&&&&&&\multicolumn{1}{c|}{}\\
\end{tabular}
}
  \caption{${\rm Ap}_1(\x ,\y ,\z )$ when $(1+\sqrt{2})t<s<(2+\sqrt{3})t$}
  \label{tb:223-1-1}
\end{table} 

From Table \ref{tb:223-1-1}, there are four candidates to take the largest value: 
\begin{align*}
&\ctext{$A$}:(s-t-1,2 s-t-1),\quad \ctext{$B$}:(s-1,2 s+t-1),\\
&\ctext{$C$}:(s+t-1,2 s-1),\quad \ctext{$D$}:(2 s+t-1,t-1)\,. 
\end{align*}
Since $s>(1+\sqrt{2})\,t$, the value at \ctext{$A$} is bigger than those at \ctext{$B$} and \ctext{$C$}.  
Because $(2 s+t-1)\x+(t-1)\y-\bigl((s-t-1)\x+(2 s-t-1)\y\bigr)=s^4-4 s^3 t+3 s^2 t^2-4 s t^3-2 t^4$ and the only real root of $x^4-4 x^3+3 x^2-4 x-2=0$ is $3.5163\dots$. 
Therefore, if $2.4142\,t=(1+\sqrt{2})\,t<s<3.5163\,t$ then the value at \ctext{$A$} is the largest and 
$$
g_1(\x,\y,\z)=(s-t-1)\x+(2 s-t-1)\y-\z\,. 
$$ 
If $3.5163\,t<s<(2+\sqrt{3})\,t=3.7320\,t$ (for example, $(s,t)=(18,5),(29,8)$) then the value at \ctext{$D$} is the largest and  
$$
g_1(\x,\y,\z)=(2 s+t-1)\x+(t-1)\y-\z\,. 
$$ 
For example, for $(s,t)=(18,5)$, $g_1(\x,\y,\z)=197871$. For $(s,t)=(29,8)$, $g_1(\x,\y,\z)=1360164$ as desired.

\section{Final comments}

When $r\ge 6$, we can also obtain the Frobenius numbers of the triple from $x^2+y^2=z^r$. However, we need to discuss the cases for each concrete value of $r$.


\end{document}